\documentclass[12pt]{article}

\usepackage{amsmath}
\usepackage{amsfonts}
\DeclareMathSymbol{\twoheadrightarrow}{\mathrel}{AMSa}{"10}
\usepackage{amssymb}
\usepackage{mathabx}
\usepackage{amsthm}     
\usepackage{url}
\title{Saturated Free Algebras and Almost Indiscernible Theories
\thanks{ 
words and phrases: Free algebras, indiscernible sets, forking 
independence, free modules}
}
\author{
\begin{tabular}{l@{\hspace{.5in}}l}
	\parbox[t]{2.5in}{	
	     T.\ G.\ Kucera\\
         Department of Mathematics\\
         University of Manitoba\\
         Winnipeg, Manitoba\\
         R3T 2N2 CANADA\\
         {\small\url{thomas.kucera@umanitoba.ca}}
         }
	&   
	\parbox[t]{2.5in}{	
         A.\ Pillay\\
         Department of Mathematics
         University of Notre Dame\\
		 Notre Dame, IN 46556 \\
		 USA\\
		 {\small\url{Anand.Pillay.3@nd.edu}}
         }
\end{tabular}\\
		}
		
\date{\today}

\newcommand{\vsp}[1]{\vspace{#1\baselineskip}}

\newcommand{\nlabel}[1]{\label{#1}}

\newtheorem{thm}{Theorem}[section]
\newtheorem{prop}[thm]{Proposition}
\newtheorem{lemma}[thm]{Lemma}
\newtheorem{cor}[thm]{Corollary}
\newtheorem{fac}[thm]{Fact}
\newtheorem{defn}[thm]{Definition}
\newtheorem{context}[thm]{Context}

\newtheorem{Rem}[thm]{Remark}
\newenvironment{rem}
   {\begin{Rem}\upshape}
   {\qed\end{Rem}}
   
\newtheorem*{Remnc}{Remark}  
\newenvironment{remnc}
   {\begin{Remnc}\upshape}
   {\qed\end{Remnc}}
   
\newtheorem{Example}[thm]{Example}
\newenvironment{exm}
   {\begin{Example}\upshape}
   {\qed\end{Example}}

\newtheorem*{Nota}{Notation}  
\newenvironment{notation}
   {\begin{Nota}\upshape}
   {\qed\end{Nota}}

\newtheorem{Query}[thm]{Question}
\newenvironment{question}
   {\begin{Query}\upshape}
   {\qed\end{Query}}

\newcounter{claimno}
\newcommand{\claimnum}
   {\noindent\stepcounter{claimno}{\emph{Claim}
\arabic{claimno}:}\quad}

\renewenvironment{proof}{\noindent\textbf{Proof:}}{\qed}   
\renewcommand{\qed}{\hspace*{\fill}\rule{1ex}{2ex}}   

\newcounter{thmpt}
\renewcommand{\thethmpt}{\textup{(\alph{thmpt})}}

\newenvironment{thmparts}
   {\begin{list}
      {\thethmpt}
      {\usecounter{thmpt}
       \leftmargin2.3em
       \labelwidth1.3em
       \labelsep0.3em
       \itemindent0em
       \listparindent0.8em
       \topsep0mm
      }
    }
   {\end{list}}

\newcommand {\setof}[2]{\left\lbrace\,#1\,\colon\,#2\,\right\rbrace}

\newcommand {\Set}[1]{\left\lbrace\,#1\,\right\rbrace}

\newcommand {\pair}[2]{\left\langle #1,\,#2\right\rangle}

\newcommand {\card}[1]{\left|#1\right|}
   
\newcommand {\gen}[1]
    {\left\langle\!\left\langle#1\right\rangle\!\right\rangle}
   
\newcommand {\tuple}[1]
    {\left\langle #1\right\rangle}
   
\newcommand {\tupleidx}[2]
    {\left\langle\,#1\,\colon\,#2\,\right\rangle}
   
\newcommand{\forkindep}[1][]{%
  \mathrel{
    \mathop{
      \vcenter{
        \hbox{\oalign{\noalign{\kern-.3ex}\hfil\(\vert\)\hfil\cr
              \noalign{\kern-.7ex}
              \(\smile\)\cr\noalign{\kern-.3ex}}}
      }
    }\displaylimits_{#1}
  }
}

\newcommand{\lsub}[2]{\vphantom{#2}_{#1}{#2}}

\newcommand{\ZZ}{\mathbb{Z}} 
\newcommand{\QQ}{\mathbb{Q}} 

\newcommand{\FF}{\mathbb{F}} 

\newcommand{\LL}{\mathcal{L}} 

\newcommand{\Mons}{\mathfrak{M}}    
\newcommand{\Upeq}{\mbox{\scriptsize eq}} 
\newcommand{\Monseq}{\Mons^{\Upeq}}       

\newcommand{\ola}{\overline{a}}
\newcommand{\olb}{\overline{b}}
\newcommand{\olc}{\overline{c}}
\newcommand{\old}{\overline{d}}
\newcommand{\olf}{\overline{f}}

\newcommand{\oly}{\overline{y}}

\newcommand{\olI}{\overline{I}}

\newcommand{\olka}{\overline{\kappa} }
\newcommand{\mub}{\bar{\mu}}

\newcommand{\newop}[1]{\expandafter\def\csname #1\endcsname
{\mathop{\mathrm{#1}}\nolimits}}

\newop{ran}            
\newop{dom}            
\newop{ann}            

\newop{Inv}            
\newop{supp}           

\newop{Th}             
\newop{pp}             
\newop{tp}             
\newop{stp}            
\newop{dcl}            
\newop{acl}            

\DeclareMathOperator{\rk}{R}

\newcommand{\ee}[1]{\underline{\mathrm{e}}_{#1}}  

\newcommand{\restr}{\upharpoonright}         

\newcommand{\tupop}{\textup{(}}
\newcommand{\tupcp}{\textup{)}}
\newcommand{\tupob}{\textup{[}}
\newcommand{\tupcb}{\textup{]}}

\begin{document}
\maketitle

\begin{abstract}
    We extend the concept of ``almost indiscernible theory''
    introduced by Pillay and Sklinos in \cite{MR3411166} (which was
    itself a modernization and expansion of Baldwin and Shelah
    \cite{MR726272}), to uncountable languages and uncountable
    parameter sequences.  Roughly speaking a theory \(T\) is almost
    indiscernible
    if some saturated model is in the algebraic closure of an
    indiscernible set of sequences.  We show that such a theory \(T\)
    is
    nonmultidimensional, superstable, and stable in all cardinals 
    \(\geq|T|\)\,.  We prove a structure theorem for sufficiently
    large
    \( a \)-models \(M\): Theorem \ref{AlmInd:prp210}
    which states that over a suitable
    base, \(M\) is in the \emph{algebraic closure} of an independent
    set
    of realizations of weight one types (in possibly infinitely many
    variables).  We also explore further the saturated free algebras
    of Baldwin and Shelah in both the countable and uncountable
    context.  We study in particular theories and varieties of
    \(R\)-modules, characterizing those rings \( R \) for which the 
	 free \( R \)-module on \( \card{R}^{+} \) generators is 
	 saturated (Theorem \ref{Free:Mod:ai}), and 
	 pointing out a counterexample to a conjecture from
    Pillay-Sklinos (Example \ref{Free:ModUTMR}).
\end{abstract}

\bibliographystyle{plain}

\tableofcontents


\section{Preliminaries}\nlabel{Pre:}   
This paper continues and builds on the investigations of Baldwin and
Shelah \cite{MR726272} and Pillay and Sklinos \cite{MR3411166}.  The
original context of Baldwin and Shelah was the study of 
\(\Th(F)\) where
\(F\) is the free algebra in \(\aleph_{1}\)-many generators in a
variety
(in the sense of universal algebra) in a countable signature or
language.  Their work was clarified (with some corrections) by
Pillay and Sklinos, and also extended to the more general notion of
\emph{almost indiscernible theories}, still in a countable language.

As stated in the abstract, roughly speaking a theory \(T\) is almost
indiscernible if some saturated model is in the algebraic closure of an
indiscernible set of sequences. In the same casual manner, one might 
say that a free algebra is one that is freely generated by a basis. 
``Algebraic closure of\ldots'' is a general model-theoretic analogue of 
``generated by\ldots''; ``indiscernible sequence''  is an analogue of 
``free (basis)''.  These analogies led Pillay and Sklinos to the 
generalizations of Baldwin-Shelah presented in \cite{MR3411166}. 

Similarly, the present abstract generalization of those results was 
in part motivated by our understanding of the structure of nicely 
behaved theories of modules: we knew that something like Theorem 
\ref{AlmInd:ssMod} should be true, given the correct definitions and 
proper development of the theory.

Although some difficult technical results in classical stability 
theory  in Section \ref{AlmInd:} provide the foundation for the entire 
paper, the ultimate goal is to pursue the applications to algebraic 
structures of various kinds.

Thus we have several aims in the current paper. First we consider the
case of almost indiscernible theories but generalized to uncountable
languages as well as indiscernible sets of infinite (rather than
finite) tuples.  Among interesting differences with the countable
case is that the theories will be superstable but not necessarily
totally transcendental.  The main structural result is 
Theorem \ref{AlmInd:prp210}.
We point out that the almost indiscernible (complete) theories \(T\)
of
modules (over some ring \(R\) in the usual language) are precisely the
superstable theories of modules which are \(|T|\)-stable.
We also revisit the special case of saturated free algebras, with
respect to a given variety in the sense of universal algebra.
Results from Pillay-Sklinos in the context of countable languages go
through smoothly for uncountable languages.  On the other hand in the
even more special case of the variety of \(R\)-modules, if the free
algebra is saturated then its theory is totally transcendental (not
just superstable).  We  classify the rings \(R\) such that the free
\(R\)-modules are saturated, and give an example where the
corresponding theory does not have finite Morley rank, yielding a
counterexample to a question from Pillay-Sklinos.

We only use very standard facts from Shelah's stability theory. 
Rather than tracing back to original sources, we rely primarily on 
the outline in Chapter 1 of Pillay \cite{Pillay:GST}, with occasional 
reference to Baldwin \cite{Bald:FunStb} where required.

In our work, we follow the proofs of \cite{MR3411166} as closely as 
possible, but give some further clarifications to the structure of 
the proofs of Propositions 2.9 and 2.10 therein, as necessitated by 
the move to the uncountable context.


\subsection{The definition}
We begin by extending Definition 2.1 of Pillay and Sklinos 
\cite{MR3411166}. 

\begin{defn}
    Let\nlabel{Pre:Df}
    \( T \) be a complete theory of cardinality \( \tau \)\,,
    \( \mu\le\tau \) a finite or infinite cardinal, 
    and \( \kappa>\tau \) a  cardinal.
    
    \( T \) is called
    \( (\mu,\kappa) \)-\emph{almost indiscernible} 
    if it has a saturated model 
    \( M \) of cardinality \( \kappa \) 
    which is in the 
    algebraic closure of an indiscernible set \( I \) of 
    \( \mu \)-sequences.
    
    \( T \) is \emph{almost indiscernible} if it is
    \( (\mu,\kappa) \)-almost indiscernible for some such
    \( \mu,\kappa \)\,.
\end{defn}

So ``almost indiscernible'' as defined in \cite{MR3411166} for 
countable theories is
\( (n,\aleph_{1}) \)-almost indiscernible for some finite \( n>0 \)\,.

Trivially if \( T  \) is \( (\mu,\kappa) \)-almost indiscernible then 
it is \( (\mu',\kappa) \)-almost indiscernible for any 
\( \mu' \)\,, \( \mu\le\mu'\le\tau \)\,, for given an indiscernible 
set of \( \mu \)-sequences, just extend each to a 
\( \mu' \)-sequence by repeating the first entry. As well, 
in the definition, we can 
replace sequences indexed by \( \mu \)  by sequences 
of cardinality \( \mu \)\,.

\begin{exm}
   Later\nlabel{Pre:not-tau} (Corollary \ref{AlmInd:Redk})
   we will see that if \( T \) 
   is almost indiscernible, it even has a saturated model of 
   cardinality \( \tau \) which is the algebraic closure of an 
   indiscernible set of \( \mu \)-tuples. 
   
   But we should not include the possibility of 
	\( \kappa =\tau \) in the definition. Let 
   \( \tuple{\mathcal{Q};\,\le} \) be the disjoint union of countably 
   many copies of the rational linear order, let \( \mathbf{e} \) be 
   some fixed enumeration of the rationals \( \QQ \) in order type 
   \( \omega \)\,, and let \( \ee{i} \) be the copy of this tuple on 
   the \( i \)-th copy of \( \QQ \) in \( \mathcal{Q} \)\,.  Then 
   clearly \( \tuple{\ee{i}:\,i<\omega} \) is an indiscernible set of 
   \( \omega \)-tuples in \( \mathcal{Q} \) whose algebraic closure 
   (in fact union) is all of \( \mathcal{Q} \)\,. 
   
   Of course the theory of this structure cannot have any of the other
   properties of an almost indiscernible theory, expounded later.
\end{exm} 

\begin{exm}
   The\nlabel{Pre:BadCh}
   definition does not require ``best possible choices''. Let 
   \(  \FF \) be an infinite field of cardinality \( \tau \)\,, 
   \( \LL \) the usual language for vector spaces over 
   \( \FF \)\,, and  \( T \) the theory of non-zero vector 
   spaces over \( \FF \)\,.  
	Let \( (\ee{\alpha})_{\alpha<\tau^{+}} \) 
   enumerate a basis of the \( \FF \)-vector space \( \mathcal{V} \)
   of dimension
   \( \tau^{+} \)\,. Then clearly 
   \( \setof{\ee{\alpha}}{\alpha<\tau^{+}} \) is an indiscernible set 
   generating \( \mathcal{V} \) and so \( T \) is 
   \( (1,\tau^{+}) \)-almost indiscernible.  
   
   But for any cardinal \( 1<\mu\le\tau \)\,, we could just as well
   take
   \( \ee{0} \) to be a \( \mu \)-tuple enumerating a basis of the 
   \( \FF \)-vector space of dimension \( \mu \)\,, and extend it to 
   a sequence \( (\ee{\alpha})_{\alpha<\tau^{+}} \) whose range is 
   again a basis of \( \mathcal{V} \)\,, exhibiting \( T \) as a 
   \( (\mu,\tau^{+}) \)-almost indiscernible theory.
   
   Alternatively, we could take \( \ee{0} \) to be an enumeration 
   of \( \FF \)\,, that is, of the one-dimensional vector space, and 
   then let \( (\ee{\alpha})_{\alpha<\tau^{+}} \)  be an enumeration 
   of  \( \tau^{+} \) direct-sum independent subspaces of 
   \( \mathcal{V} \)\,, exhibiting \( T \) as a 
   \( (\tau,\tau^{+}) \)-almost indiscernible theory, with a lot more 
   information than is really required.
  
   One of the goals of the structure theory we develop is to recover 
   some of the fine detail that might be lost by redundancy in the 
   indiscernible set of \( I \)-sequences.
\end{exm}

\begin{exm}
   We\nlabel{Pre:nonAI} point out briefly that both the condition of
   almost indiscerniblity, as well as the consequences in Theorem
   \ref{AlmInd:prp210} below, concern exceptional behaviour, even for
   uncountably categorical theories \(T\)\,.  Of course when \(T\) is
   almost strongly minimal, that is, any model is in the algebraic
   closure of a fixed strongly minimal set \(D\) (without parameters
   say), then we do have almost indiscernibility: any model \(M\) is
   in the algebraic closure of the indiscernible set consisting of a
   maximal independent set of realizations in \(M\) of the generic
   type of \(D\)\,.  Even when \(T\) is not almost strongly minimal,
   such as \(\Th(\ZZ/4\ZZ)^{(\omega)}\)\,, the conditions of almost
   indiscernibility may still hold.  Let \(1\) be the generator of a
   copy of \(\ZZ/4\ZZ\)\,: then \(\tp(1)\) has Morley rank \(2\) but
   weight \(1\)\,, and of course any model is generated as a
   \(\ZZ/4\ZZ\)-module by an independent set of realizations of this
   type.
      
   However consider the theory \(T\) of the structure consisting of
   two sorts \(X\) and \(V\) with surjective \(\pi:X\to V\)\,, and
   where
   \(V\) has the structure of an infinite-dimensional vector space
   over \(\FF_{2}\) say, and each fibre \(X_{a}\) is (uniformly
   definably) a principal homogeneous space for \((V,+)\)\,.  Clearly 
   \( T \) can be axiomatized in a two-sorted language with symbols 
   \( \pi \) of sort \( X\rightarrow V \)\,, 
   \( + \) of sort \( V\times V\rightarrow V \)\,, 
   and \( \pair{\ }{{}}\) of sort \( V\times X\rightarrow X\)\,.  
   Then \(T\) is \(\aleph_{1}\)-categorical, but
   is {\em not} almost indiscernible.  \(V\) is strongly minimal and
   \(X\) has Morley rank \(2\)\,, degree \(1\)\,.  Let \(a\) be a
   generic point of \(V\) and \(b\in X\) be such that 
	\(\pi(b) = a\)\,.
   Then \(\tp(b)\) has weight \(1\)\,.  If \(\setof{b_{i}}{i\in I}\)
   is
   a maximal independent set of realizations of \(\tp(b)\) in a model
   \(M\)\,, then necessarily \( \setof{\pi(b_{i})}{i\in I} \) is a
   linearly independent set in \( V \)\,.  Furthermore, 
   \( \acl(V)=  V\)\,, so for distinct non-zero
   \( u,\,v\in V \)\,, no part of the fibre over
   \( u+v \) is algebraic over the (union of the) 
	fibres over \( u \) and \( v \)\,.
   So on the one
   hand any maximal indiscernible set \( I  \) (of tuples) in a 
   model \( M \)
   cannot intersect all the fibres of \( \pi \)\,, and on the other 
   hand any fibre that does not intersect \( I \) is not in the 
   algebraic closure of \( I \)\,.
\end{exm}


\subsection{Basic facts}

We need some translations of the basic facts about stability theory 
enunciated in Section 1 
of \cite{MR3411166}. Our theories will be superstable, not necessarily
totally transcendental,
and so may not have prime models.  As is usual, we will abbreviate 
``superstable'' as ``ss'' and ``totally transcendental'' as ``tt''. 
We remind the reader that for \emph{uncountable} languages we have to 
characterize tt theories as those where every formula (every type) 
has ordinal-valued Morley rank. For countable theories only, this is 
equivalent to \( \omega \)-stability.
But each theory \( T   \) that
we consider will nonetheless be stable in \( \tau=\card{T} \)\,, and 
will have a model 
\( M_{\omega} \) which is \( \omega \)-saturated,  is an 
\( a \)-model, and such that every stationary type is non-orthogonal 
to a type over \( M_{\omega} \)\,.  
 Furthermore, 
\emph{every} model we care about is an elementary extension of 
\( M_{\omega} \)\,, so there are still very strong parallels to 
\cite{MR3411166}. The principal difficulty lies not so much in the 
movement to merely superstable theories, but in allowing infinite 
sequences as the elements of the indiscernible sets.

\vsp{1}
For the remainder of this section and in Section 2 (unless explicitly 
stated otherwise), \( T \) is a 
\textbf{complete superstable theory}  in a
language of 
cardinality \( \tau \)
and \( \Mons \) is a sufficiently large saturated model of 
\( T \) (a universe), with \( \Monseq \) being the associated 
``imaginary'' universe. We  work in \( \Monseq \): 
every element, set, sequence, model that we 
consider is a ``small'' thing in \( \Monseq \)\,, that is, of 
cardinality strictly less than the cardinality of \( \Mons \)\,. 
By ``algebraic (or definable) closure'' we 
always mean ``in the sense of \( \Monseq \)''.  However, 
in reading the 
algebraic examples, it is always helpful to think of things taking 
place in the ``home'' sort.

\begin{fac}\textup{\cite[4.1.1, 4.1.2, 4.2.1]{Pillay:GST}}
   There is\nlabel{Pre:SS} a cardinal \( \lambda(T)\le 2^{\tau} \) 
   such that \( T \) is stable in \( \kappa \) if{}f 
   \( \kappa\ge\lambda(T) \)\,. \( T \) has a saturated model 
   in every cardinal \( \kappa\ge\lambda(T) \)\,. Since \( T \) is 
   superstable, \( \kappa(T)=\aleph_{0} \)\,, that is, every type 
   \tupop in finitely many variables\tupcp\ 
   does not fork over some finite set.
\end{fac}

Recall that if \(A, B, C\) are sets of parameters (or tuples),  
then \(B\)
\emph{dominates} \(C\) \emph{over} \(A\)
if whenever \(D\) is independent from
\(B\) over \(A\) then \(D\) is independent from \(C\) over \(A\)\,.

\begin{defn}{\quad}\nlabel{Pre:aMod}
   \begin{thmparts}
      \item The \emph{strong type} of \( a \) over \( A \)\,, 
      \( \stp(a/A) \)\,, is the type of \( a \) over \( \acl(A) \)
      \tupop for emphasis, in \( \Monseq \)\tupcp.

      \item By an \( a \)-model of \(T\) we mean a model \(M\) of
      \(T\)
      such that any strong type over any finite subset of \(M\) is
      realized in \(M\)\,.
   
      \item A type \(p(x)\in S(A)\) is said to be \( a \)-isolated if
      there is a finite subset \(B\) of \(A\) and a strong type
      \( q(x) \) over \(B\) which implies \(p(x)\)\,.
   \end{thmparts}
\end{defn}

\begin{fac}\nlabel{Pre:aPrm}\textup{\cite[4.2.4, 4.3.4]{Pillay:GST}}
   \begin{thmparts}
      \item For
      any set of parameters \(A\) there is an \( a \)-prime
      model over \(A\)\,, that is, an \( a \)-model \(M\) containing
      \(A\)
      such that for any \( a \)-model \(N\) containing \(A\) there
      is an
      elementary embedding over \(A\) of \(M\) into \(N\)\,.  \(M\)
      has
      the property that for all tuples \(b\) from \(M\)\,,
      \(tp(b/A)\) is
      \( a \)-isolated.
   
      \item Suppose
      \(M_{0}\) is an \( a \)-model, and \(A\) is any set
      of parameters and \(b\) any tuple.  Then \(tp(b/M_{0}A)\) is
      \( a \)-isolated iff \(A\) dominates \(Ab\) over \(M_{0}\)\,.
   \end{thmparts}
\end{fac}

Clearly an \( a \)-model is \( \aleph_{0} \)-saturated.

\begin{defn}
 \(T\)\nlabel{Pre:nmd} is \emph{nonmultidimensional} 
 \textup{[``nmd'']} if every
stationary type \(p\) is nonorthogonal to \(\emptyset\)\,, that is,
nonorthogonal to some stationary type which does not fork over
\(\emptyset\)\,. 
\end{defn}

\begin{rem}
Definition\nlabel{Pre:nmd:r} \ref{Pre:nmd}  
is equivalent to the following: 
\begin{quote}\em
   \noindent For any \( A \) and any
   stationary type \(q(x) \) over \( A \)\,, if 
   \( \stp(A'/\emptyset) =\stp(A/\emptyset)\)\,, 
   \(A'\) is independent from \(A\) over \(\emptyset\)\,,  and
   \(q'\) is the copy  of \(q\) over \(A'\)\,, then \(q\) is
   nonorthogonal to \(q'\)\,. 
\end{quote}

In fact, at least for superstable theories, in \ref{Pre:nmd} it 
suffices to demand that every stationary type \(p\) be non-orthogonal 
to a type over a fixed \( a \)-model (as mentioned in the first 
paragraph of this Section), cf.\ Baldwin 
\cite[XV: Theorem 1.7]{Bald:FunStb}.
\end{rem}

We need the following  basic result:

\begin{prop}
   Let\nlabel{Pre:ext-a-m}
   \(T\) be a superstable nonmultidimensional
theory \tupop of any cardinality\tupcp. Then any elementary extension
of an
\( a \)-model is an \( a \)-model.
\end{prop}

This proposition is folklore and there are various routes to it.  For
example it follows directly from Shelah's `three model lemma', and
also
follows from Propositions 3.2 and 3.6 of Chapter 7 of
\cite{Pillay:GST}.  We will give a quick independent proof, starting
with a suitable \(3\)-model lemma.

\begin{lemma} 
Suppose\nlabel{Pre:3mod} \(T\) is superstable nonmultidimensional,
and  \(M_{0}\prec M\precneq N\) where \(M_{0}\) is an \( a \)-model.
Then there is \(c\in N\setminus M\) such that \(tp(c/M)\) is regular 
and does not fork over \(M_{0}\)\,.
\end{lemma}

\begin{proof} Choose \(b\in N\setminus M\) such that
\(\rk^{\infty}(tp(b/M)) = \alpha\) is minimized.  Let
\(\varphi(x,a)\) with
\(a\in M\) be a formula in \(tp(b/M)\) of \(\infty\)-rank
\(\alpha\)\,. In
particular \(tp(b/M)\) does not fork over \(a\)\,. As \(M_{0}\) is an
\(a\)-model we can choose \(a'\in M_{0}\) such that 
\(\stp(a') = \stp(a)\)
and \(a'\) is independent from \(a\) over \(\emptyset\)\,. By Remark 
\ref{Pre:nmd:r},
there is a type \(q(x)\) over \(M_{0}\) which contains the formula
\(\varphi(x,a')\) and is nonorthogonal to \(tp(b/M)\)\,. So there is
\(M'\supseteq M\) with \(b\) independent from \(M'\) over \(M\) and
\(c'\)
realizing \(q|M'\) such that \(b\) forks with \(c'\) over \(M'\)\,.

A standard
argument yields \(c\in N\setminus M\) such that 
\(\models \varphi(c,a')\)\,:
There is a formula \( \chi(x,y,z) \)  over \( M \)
and \(d\in M'\) such that \(\models \chi(b,c',d)\) 
witnesses the forking of
\(b\) with \(c'\) over \(M'\)\,, that is,
\(\chi(x,c'',d'')\cup tp(b/M)\) forks
over \(M\) for any \(c'',d''\)\,.  Now \(\exists y(\chi(x,y,d)\wedge
\varphi(x,a'))\) is in \(tp(b/M')\) so for some \(d'\in M\)\,,
\(\exists
y(\chi(x,y,d')\wedge \varphi(x,a'))\) is in \(tp(b/M)\)\,. Let \(c\in
N\) be
such that \(\models \chi(b,c,d')\wedge\varphi(c,a')\)\,.
As \(b\) forks with
\(c\) over \(M\)\,, \(c\in N\setminus M\)\,.

So as \(R^{\infty}(\varphi(x,a')) =\alpha\)\,, 
by the minimal choice of \(\alpha\)\,, 
\(R^{\infty}(tp(c/M) = \alpha\)\,.
Then as \(a'\in M_{0}\)\,, \(tp(c/M)\) does not fork over \(M_{0}\)\,,
as required. Furthermore, \(tp(c/M)\) is regular by the choice of 
\( \alpha \) minimal, as in the argument of
\cite[Lemma 4.5.6]{Pillay:GST}.
\end{proof} 

\vsp{1}
\begin{proof}
   (of Proposition \ref{Pre:ext-a-m})
   
\(T\) is assumed to be superstable  nonmultidimensional. 
Let \(M \) be an \( a \)-model, and assume
\(M \prec N\)\,. We want to prove that \(N\) is an \(a\)-model. 
Let \(N'\) be the \(a\)-prime model over \(N\) given by Fact
\ref{Pre:aPrm}(a).  It will be enough to show that \(N = N'\)\,.

Suppose not. Then by Lemma \ref{Pre:3mod}, there is 
\(c\in N'\setminus N\) which is independent from \(N\) over \(M\)\,.
But
\(tp(c/N)\) is \( a \)-isolated, so by Fact \ref{Pre:aPrm}(b), 
\(N\) dominates \(c\)
over \(M \) which is a contradiction. 
\end{proof}

\begin{cor}
	\tupop\( T \) superstable nonmultidimensional\tupcp\ 
   If\nlabel{Pre:PrmMin}
   \( M \) is an \( a \)-model and \( N \) is \( a \)-prime 
   over \( M \cup A \)\,, then \( N \) is prime and minimal over 
   \( M \cup A \)\,.
\end{cor}

\begin{proof}
	 It is immediate by \ref{Pre:ext-a-m} that \( N \) is prime; and 
	 by the proof just given, it is minimal over 
	 \( M \cup A \)\,.
\end{proof}


\vsp{1}

Now let \( \mathcal{A}_{0} \) be the \( a \)-prime model of \( T \)  
 over  \( \emptyset \)\,. 

Let \( (p_{i})_{i\in I} \) be a list, up to non-orthogonality, of all 
the regular types over \( \mathcal{A}_{0}  \) (and hence, up to 
non-orthogonality, all the regular types of \( T \)). Since \( T \) 
is superstable, for each \( i\in I \) there is finite 
\( a_{i} \in \mathcal{A}_{0} \) such that \( p_{i} \) 
is definable over
\( a_{i} \)\,. [Since \( T \) is superstable,  we can 
choose finite \( b\in \mathcal{A}_{0} \) so that \( p_{i} \) does not 
fork over \( b \)\,; then since \( \mathcal{A}_{0} \) is an 
\( a \)-model, we can find \( c\in\mathcal{A}_{0} \) realizing
the restriction of \( p_{i} \) to \( b \) and in the correct strong 
type. Then  \( p_{i} \) is definable over \( a_{i}=bc \)\,.]
We let \( \hat{p}_{i} \) be the restriction of
\( p_{i} \) to \( \acl(a_{i}) \)\,. 

\begin{remnc}
   Note that in general \( \card{I}\le 2^{\card{T}} \)\,, but we will 
   see that in the 
   context we develop, as \( T \) has an a-prime model of 
   cardinality \( \tau=\card{T} \) and \( T \) is \( \tau \)-stable, 
   in fact \( \card{I}\le\tau \)\,.
\end{remnc}

We need a slight reformulation of  \cite[Fact 1.3]{MR3411166}.

\begin{lemma}\nlabel{Pre:13+}
    Let \( \mathcal{A}_{0}\prec M \prec \Mons\)\,. 
    For each \( i\in I \) let 
    \( J_{i}  \) be a maximal independent set of realizations of 
    \( p_{i} \) in \( M \)\,. Then \( M \) is a-prime, prime, and
    minimal 
    over \( \mathcal{A}_{0}\cup\bigcup_{i\in I} J_{i} \)\,.
\end{lemma}


\section{Almost indiscernible theories}\nlabel{AlmInd:}   

\begin{context} 
	 {\upshape Unless explicitly stated otherwise, for Section 2:}
    \begin{thmparts}{}\nlabel{AlmInd:Ctx}
        
        \item  \( T \) is a \( (\mu,\kappa) \)-\emph{almost
        indiscernible} theory, \( \card{T}=\tau \)\,, 
        \( \mu\le\tau<\kappa \)\,,
        with universe  
        \( \Mons \) of some regular cardinality 
		  \( \olka\gg\kappa \)\,.

        \item  \( M \) is a saturated model as in the definition: 
        \( \card{M}=\kappa \)\,, \( I \) is an indiscernible 
        \emph{set}  of 
        \( \mu \)-sequences in \( M \)\,, and \( M \) is in the 
        algebraic closure of the \tupop union of\tupcp\ \( I \)\,.  
        
        Since 
        \( \mu\le\tau<\kappa \)\,, necessarily 
        \( \card{I}=\kappa \)\,, so we can write \( I \) as a 
        \( \kappa \)-sequence 
        \( \tuple{\ee{\alpha}:\alpha<\kappa} \)\,, and when
        necessary the \( \mu \)-sequence \( \ee{\alpha} \) is indexed
        as 
        \( \tuple{\ee{\alpha,i}:i<\mu} \)\,.
        
        \item Extend \( I \) to an indiscernible `set' 
        \( \olI=\tuple{\ee{\alpha}:\alpha<\olka} \) in \( \Mons \)\,.

        For each infinite \emph{ordinal} \( \lambda\le\olka \)\,, let 
        \( I_{\lambda}=\tuple{\ee{\alpha}:\alpha<\lambda} \) and set 
        \( M_{\lambda}=\acl(I_{\lambda}) \) in \( \Mons \)\,.  Note 
		  the extension of this sequence to the size of the universe.
    \end{thmparts}
\end{context}
        
In particular, \( M_{\kappa}=M \) is an elementary substructure of 
\( \Mons \)\,, but the status of all the other \( M_{\lambda} \) 
remains to be resolved.


\subsection{Basic Facts}

Recall that for an infinite cardinal \( \nu \)\,,
 \( M \) is \( F^{a}_{\nu} \)-saturated 
 if every strong type over any subset of \( M \) of 
 cardinality less than \( \nu \) is realized in 
 \( M \)\,.

\begin{thm}{\quad}
     \( \lambda \)\nlabel{AlmInd:Ch}
	  denotes an infinite  \emph{ordinal}. 
    \begin{thmparts}
        \item 
		  \( \lambda \geq \tau \)\label{ca}  implies 
		  \( \card{M_{\lambda}} = \lambda \) and 
		  \( \lambda < \tau \) implies
		  \( \mu\card{\lambda} \leq \card{M_{\lambda}}\leq \tau \)\,.

        \item  
		  For\label{cb} all \(\lambda \leq \olka \)\,, 
		  \(M_{\lambda} \preceq \Mons\ \)\,. 
		  
		  \tupop And so 
		  \( \tuple{M_{\lambda}}_{\omega\le\lambda\le\olka} \) 
		  is an elementary chain.\tupcp
    
        \item  
		  For all \label{cc}
		  \( \lambda \leq \olka \), \( M_{\lambda} \) is
		  \( \card{\lambda} \)-saturated.  
		  
        \item  In particular \label{cd} \( M_{\olka} \) is saturated 
		  of  cardinality \( \olka \)\,, so without loss of generality 
		  equal to \( \Mons \)\,.  
        
        \item For all\label{ce} \( \lambda \leq \olka \)\,, 
		  \(M_{\lambda}\)  is \( F^{a}_{\card{\lambda}} \)-saturated.
        
		  \item In\label{cf} particular\nlabel{AlmInd:FaMod}
		  all \(M_{\lambda}\) have the 
		  property that all strong types over finite sets are realized.
		  
		  Then it will follow from Theorem \ref{AlmInd:Tss} to come, 
		  that since \( T \) is superstable, 
		   all the \(M_{\lambda}\) are \(a\)-models.
    \end{thmparts}
\end{thm}

\begin{proof}
\begin{thmparts}
	 \item follows by simple counting.
	 
	 \item  By definition \( M = M_{\kappa} \) is an elementary 
	 substructure of \( \Mons \)\,.  
	 First assume  \( \lambda < \kappa \)\,.  Then an easy Tarski-Vaught
	 argument which we now describe yields that  
	 \( M_{\lambda} \preceq M_{\kappa} \)\,: 
	 Let \( \varphi(x) \) be a formula with parameters from
	 \( M_{\lambda} \) with a solution \( d\in M_{\kappa} \)\,. 
	 Let \(J\) be a finite
	 subset of \( \lambda \) and \( \alpha_{1},\ldots,\alpha_{n} \) 
	 distinct elements 
	 of \( \kappa\setminus\lambda \), and such that the parameters in
	 \( \varphi(x) \) are from 
	 \( E = \acl(\setof{\ee{\beta}}{\beta\in J}) \) and  
	 \( d\in \acl(E\cup\setof{\ee{\alpha_{i}}}{i=1,\ldots,n})\). 
	 Let \(\gamma_{1},\ldots,\gamma_{n}\) be distinct elements of 
	 \(\lambda\setminus J\)\,. 
	 Then by indiscernibility, 
	 \( \tp(\ee{\gamma_{1}},...,\ee{\gamma_{n}}/E) =
	 \tp(\ee{\alpha_{1}},\ldots,\ee{\alpha_{n}}/E) \)
	 (where these types are computed in  \( M_{\kappa}\)).
	 Thus  we can find
	 \( d'\in \acl(E, \ee{\gamma_{1}},\ldots,\ee{\gamma_{n}}) 
	 \subseteq M_{\lambda} \)
	 such that \( M_{\kappa}\models \varphi(d') \)\,.

	 On the other hand,  assume that  \(\lambda > \kappa\)\,.  We will 
	 also use Tarski-Vaught: 
	 Let \( \varphi(x) \) be a formula over \( M_{\lambda} \) such that 
	 \( \Mons\models \exists x \varphi(x) \)\,.  
	 Let \( J \) be a finite subset of \( \lambda \) such that the 
	 parameters of \( \varphi \) 
	 are in \(E =\acl(\setof{\ee{\beta}}{\beta\in J}) \)\,.  
	 Choose a subset \( J' \) of
	 \( \kappa \) of the same cardinality as \( J \)\,, 
	 so by indiscernibility \( J \) and \( J' \) have the same type 
	 (under any enumeration) in \( \Mons \)\,. 
	 Let \( E' = \acl(J') \)\,, and let \( \varphi'(x)\) be the image of 
	 \( \varphi(x) \) under a partial elementary map \( f \) taking 
	 \( J \) to \( J' \) and \( E \) to \( E' \)\,. 
	 As \( M_{\kappa}\prec \Mons \)\,, 
	 \( M_{\kappa}\models \exists x \varphi'(x) \)\,.
	 So let \( d'\in M_{\kappa} \) be such that
	 \( M_{\kappa} \models \varphi(d') \)\,.
	 Again \( d' \) will be in \( \acl(J', K') \) 
	 where \( K' \) is a finite subset of \( \kappa \) disjoint from 
	 \( J' \)\,. Let \( K \) be a finite subset of \( \lambda \)
	 of the same cardinality as \(K'\) and disjoint from \( J \)\,.  
	 Then again indiscernibility implies that  the partial elementary 
	 map \( f \) (in the sense of \( \Mons \)) extends to a 
	 partial elementary map \( g \) taking \( K \)
	 to \( K' \) and \( \acl(E,K) \) to \( \acl(E',K') \)\,.  
	 Let \( g(d) = d' \)\,, 
	 so \( \Mons \models \varphi(d)\) and \( d\in M_{\lambda} \)\,.  

	 \item Note that this part is an analogue of Morley's theorem 
	 (that for a countable complete theory \(T\)\,, if for some 
	 uncountable cardinal \( \kappa \) all models of \( T \) of 
	 cardinality  \( \kappa \) are saturated, 
	 then all uncountable models of \( T \) are saturated), and our 
	 proof of 
	 the harder case (where \( \lambda >\kappa \)) will be closely 
	 related to the  proof of Morley's theorem as given in
	 \cite[Theorem 5.33]{Pillay:ModTh}.
	
	 The first (and easy) case is when \( \lambda < \kappa \)\,.
	 Then the proof of part \ref{cb} (in the case 
	 \( \lambda < \kappa \))  adapts. Namely in this
	 case we have a type \( p(x) \) over a set of parameters 
	 \( E\subset M_{\lambda} \) of cardinality \( <\card{\lambda} \)\,, 
	 and now choose \( J\subset\lambda \) of cardinality 
	 \( <\card{\lambda} \)
	 with \( E\subseteq \acl(J)\)  (without loss \( E = \acl(J) \)).  
	 \( p \) is realized in \( M_{\kappa} \) by some \( d \) in
	 the algebraic closure of \( E \) together with finitely many
	 \( \ee{\alpha} \) with \( \alpha\in \kappa \setminus J \)\,.  
	 Then as \( \card{J}<\card{\lambda} \)\,, we can again replace 
	 these \( \ee{\alpha} \) by some \( \ee{\gamma} \) in
	 \( M_{\lambda} \) and realize \( p \) in \( M_{\lambda} \)\,.

	 The harder case occurs when \(\lambda > \kappa\)\,.   

	 There is no harm (for notational simplicity) in assuming that
	 \( \lambda \) is a cardinal. We  suppose that \( M_{\lambda} \) 
	 is  not saturated, and aim for a contradiction. Then there is a 
	 subset \( A \) of  \( M_{\lambda} \) with
	 \( \card{A} <\lambda \) and a complete type \( p(x) \) over
	 \( A \) which is not realized in \( M_{\lambda} \)\,. 
	 We may assume that \( A \) has cardinality at least \( \tau \)\,. 
	 Let \( J\subseteq I_{\lambda} \) be of
	 cardinality \( \card{A} \) such that \( A\subseteq \acl(J) \)\,, 
	 and let \( I \) be a countable subset of \( I_{\lambda} \) 
	 disjoint from \( J \)\,.   Extending \( p \)
	 to a complete type over \( \acl(J) \) we may assume that 
	 \( A = \acl(J) \)\,.
	 Note that \( I \) is indiscernible over \( A \)\,.
	 Note also that, by part \ref{cb}
	 \( \acl(A,I) \) is an elementary substructure \( N \) of 
	 \( \Mons \) (and of \( M_{\lambda} \)), and \( p(x) \) is not
	 realized in  \( N \)\,.  In particular
	 for any consistent formula \( \varphi(x) \) over 
	 \( A\cup I \)\,, we can pick a
	 formula \( \psi_{\varphi}(x)\in p(x) \) such that 
	 \(\models \exists x(\varphi(x)\wedge \neg\psi_{\varphi}(x)) \)\,. 
	 (Here ``\( \models \)'' means, equivalently, in
	 \( \Mons \) or in  \(M_{\lambda}\) or in \(N\)\,.)

	 We now construct a subset \( J' \) of \( J \) of cardinality
	 \( \tau \)\,, and 
	 \( A'= \acl(J')\subseteq A \) such that 
	 for \( p'(x) = p|A' \) we have:
	 
    \begin{quotation}
	 	  \noindent(*) For each consistent formula \( \varphi(x) \) 
		  over \( A' \cup I \)\,, there is
	 	  a formula \( \psi(x)\in p'(x) \) such that 
		  \(\models \exists x(\varphi(x)\wedge\neg\psi(x)) \)\,.
    \end{quotation}

	 We do this by a routine union of chain argument. 
	 We define a sequence of pairs \( \pair{J_{i}}{A_{i}} \)\,,
	 \( J_{i}\subseteq J \)\,, \( \card{J_{i}}=\tau \)\,, and 
	 \( A_{i}= \acl(J_{i})\subseteq A \) by recursion on 
	 \( n<\omega \)\,.

	 Let \( J_{0} \)
	 be any subset of \( J \) of cardinality \( \tau \)\,, and 
	 \( A_{0} = \acl(J_{0})\subseteq A \)\,. 
	 
	 Given \( \pair{J_{i}}{A_{i}} \)\,,
	 for each consistent formula \( \varphi(x)\) over
	 \( A_{i}\cup I \)\,, we have \( \psi_{\varphi}(x)\in p(x) \) 
	 such that 
	 \( 
	 \models\exists x(\varphi(x)\wedge \neg \psi_{\varphi}(x)) 
	 \)\,. 
	 Note that there are at
	 most \( \tau \) such formulas \( \psi_{\varphi}(x) \)\,.
	 Add the parameters (from \( A \)) of the formulas 
	 \( \psi_{\varphi}(x) \) to \( A_{i} \) to obtain
	 \(  A_{i}' \) (which still has cardinality \( \tau \)). 
	 Extend \( J_{i} \) to \( J_{i+1}\subseteq J \) of cardinality 
	 \( \tau \) such that \( A_{i}'\subseteq \acl(J_{i+1}) \)\,, 
	 and set \( A_{i+1}=\acl(J_{i+1}) \)\,.  
	 Set \( A' = \bigcup_{n} A_{n}\) and 
	 \( J' = \bigcup_{n} J_{n}\)\,.   
	 Then set \( p'(x) \) to be the restriction of \( p \) to
	 \( A' \)\,.  
	 So we have obtained (*). 

	 Now let \( I' \) be a subset of \( I_{\lambda} \) of cardinality 
	 \( \kappa \) which is disjoint from \( J' \)\,. 
	 Then \( \acl(J',I') \) is an elementary
	 substructure \( M' \) of \( \Mons \) which is isomorphic to
	 \( M_{\kappa} \) (as \( \card{J'\cup I'} = \kappa \)).  
	 By hypothesis \( M' \) is
	 \( \kappa \)-saturated, and \( A' \) has cardinality 
	 \(\tau < \kappa\), so the type \( p'(x) \) is realised in 
	 \( M' \)\,, by some \(d'\)\,.
	 But \( d'\in acl(A',I') \) so its type over \( A'\cup I' \) is 
	 isolated by a consistent algebraic formula \( \varphi(x) \)
	 over \( A'\cup I' \)\,.  We will
	 exhibit the parameters from \( I' \) by writing \( \varphi \) as
	 \( \varphi(x,b_{1},\ldots,b_{n})\)  where  
	 \( b_{i} \) is a finite tuple from some
	 \( \ee{i} \in I' \)\,. 
	 Now  
	 \( M\models 
	 \forall x (\varphi(x,b_{1},..,b_{n})  \to \psi(x))\) 
	 for all \( \psi(x)\in p'(x) \)\,, as 
	 \( \varphi(x,b_{1},\ldots,b_{n}) \) isolates
	 \( \tp(d/A'\cup I') \) and \( d' \) realizes \( p \)\,. 

	 But  
	 \( \tp(b_{1},\ldots,b_{n}/A') = \tp(c_{1},\ldots,c_{n}/A')\)
	 for some finite tuples \( c_{i} \) from 
	 \( \ee{i}\in I \)\,.  
	 But then \( \varphi(x,c_{1},\cdots,c_{n}) \) is
	 consistent, and we have: 
	 \( \models \forall x(\varphi(x,c_{1},\ldots,c_{n})\to \psi(x)) \)  
	 for all \( \psi(x)\in p'(x) \)\,, which contradicts (*).  
	 This completes the proof of part \ref{cc}. 

	 \item Immediate.  \(M_{\olka}\) is now known to be a
	 \( \olka \)-saturated model, of cardinality \( \olka\)\,, so
	 can be assumed to be the monster model \( \Mons \)\,.

	 \item Note that a strong type over a set \( A \) is precisely a 
	 type over \( \acl^{\Upeq}(A) \)\,. 
	 We can repeat the proof of  \ref{cc} working now  with
	 types over algebraically closed sets in \( \Monseq \)\,.
	 Alternatively, at least for \( A \) the (real) algebraic closure of 
	 an infinite subset of \( \olI \)\,, \( A \) is already an elementary
	 substructure of \( \Mons \)\,,
	 so complete types over \( A \) and strong types
	 over \( A \) amount to the same thing. 
	 
	 \item is just a specialization of \ref{ce}.

\end{thmparts}    
    
\end{proof}

In particular, for cardinals \( \nu\ge\tau \)\,,
\( \tuple{M_{\lambda}:\,\nu\le\lambda<\nu^{+}} \) is a elementary 
chain of copies of the saturated model of \( T \) of cardinality 
\( \nu \)\,.

\begin{cor}
    Let\nlabel{AlmInd:Redk}
    \( T \) be a complete theory of cardinality \( \tau \)\,.
    \begin{thmparts}
       \item Then\label{AIRa}
		 \( T \) is \( (\mu,\kappa) \)-almost indiscernible
       for some \( \kappa>\tau \) if{}f \( T \) is 
       \( (\mu,\tau^{+}) \)-almost indiscernible if{}f \( T \) is 
       \( (\mu,\kappa')\)-almost indiscernible for all
       \( \kappa'>\tau \)\,.
    
       \item  In particular, under the conditions of \ref{AIRa},
		 \( M_{\tau} \) is a saturated model 
       which is the algebraic closure of 
       an indiscernible sequence 
		 \( \tuple{\ee{\alpha}:\alpha<\tau} \) 
       of \( \mu \)-tuples\,.
       
    \end{thmparts}
\end{cor}

\begin{remnc}
    We can now assume without loss of generality that 
    \( T \) is a complete theory of cardinality \( \tau \) which is
    \( (\mu,\tau^{+}) \)-almost indiscernible for some 
    \( \mu\le\tau \)\,.
\end{remnc}

\begin{thm}
    Let\nlabel{AlmInd:Tss}
    \( T \) be \( (\mu,\tau^{+}) \)-almost indiscernible. 
    Then \( T \) is stable in every cardinal
    \( \lambda\ge\tau \)\,, hence \( T \) is superstable, and in 
    particular if  \( T \) is countable and
     \( (\aleph_{0},\aleph_{1}) \)-almost indiscernible, then 
    \( T \) is tt.
\end{thm}

\begin{proof}
    Once again, we show that \( T \) is stable in all cardinals
    \( \lambda\ge\tau \)\,, by the method of the proof of
    \cite[Proposition 2.5]{MR3411166}.
    
    So let \( \lambda\ge\tau \) be a cardinal. Since 
	 \( M_{\lambda}  \) is 
    saturated of cardinality \( \lambda \)\,, it suffices to count 
    the complete types over \( M_{\lambda}  \)\,. Let 
    \( p(v) \) be a complete type over \( M_{\lambda}  \)\,.
    Then \( p \) is realized in \( M_{\lambda^{+}} \) by some element 
    \( d \)\,, which is then algebraic over 
    \( M_{\lambda}\cup
    \Set{\ee{\alpha_{1}},\,\ldots,\,\ee{\alpha_{n}}} \) 
    with 
    \( \lambda\le\alpha_{1}<\,\cdots\,<\alpha_{n}<\lambda^{+} \)\,.
    So the type of \( d \) over \( M_{\lambda}\cup
    \Set{\ee{\alpha_{1}},\,\ldots,\,\ee{\alpha_{n}}} \) is isolated
    by 
    some formula \( \theta(v,c_{1},\,\ldots,\,c_{k}) \)\,, where
    \( \theta(v,x_{1},\,\ldots,\,x_{k}) \) is a formula with 
    parameters from \( M_{\lambda} \) and 
    \( c_{1},\,\ldots,\,c_{k} \) are entries from the sequences
    \( \ee{\alpha_{1}},\,\ldots,\,\ee{\alpha_{n}} \)\,.
    There are at most \( \lambda \) many such formulas
	 \( \theta \)\,.
    By indiscernibility, the type of 
    \( \tuple{\ee{\alpha_{1}},\,\ldots,\,\ee{\alpha_{n}}} \) over
    \( M_{\lambda}  \) depends only on \( n \)\,. There are only
    \( \mu<\lambda \) choices for finite sequences 
    \( c_{1},\,\ldots,\,c_{k} \) from 
    \( \tuple{\ee{\alpha_{1}},\,\ldots,\,\ee{\alpha_{n}}} \)\,.
    Hence there are for fixed such \( \theta \)\,, no more than 
    \( \mu \) possibilities for the type of \( d \) over 
    \( M_{\lambda}\cup
    \Set{\ee{\alpha_{1}},\,\ldots,\,\ee{\alpha_{n}}} \)\,, so
    certainly 
    no more than \( \mu \) types over \( M_{\lambda} \) whose 
    realizations are determined by \( \theta \)\,. Therefore there 
    are no more than \( \mu\lambda=\lambda \) 1-types over
    \( M_{\lambda} \)\,.
\end{proof}

\begin{remnc}
    For uncountable \( T \)\,, all that follows in general for 
    superstable theories is that \( T \) is stable in every cardinal 
    \( \ge 2^{\card{T}} \)\,.  So almost indiscernible theories are 
    ``strongly''
    superstable. Any complete theory \( T \) in a (possibly 
    uncountable language) which is categorical in \( \card{T}^{+} \) 
    is ``strongly'' superstable, cf.\ the revised edition of Shelah's 
    ``Classification Theory'', the  first paragraph of the 
    proof of \cite[Theorem IX.1.15]{MR1083551}.
\end{remnc}

\vsp{1}
We need to make the use of infinitely many variables precise.
    Let
    \( \vec{v}=\tuple{v_{\beta}:\beta<\mu} \) be a sequence of 
    distinct variables. A formula \( \varphi(\vec{v}) \)
    ``in the variables \( \vec{v} \)'' 
    is some finitary formula \( \varphi \) 
    in some \tupop definite\tupcp\ finite list of 
    variables from  \( \vec{v} \)\,. This establishes a  
    correspondence between formulas and their variables,
	 and \( \mu \)-sequences of 
    elements, for the purposes of definitions such as the following:

\begin{defn}    
    Set\nlabel{AlmInd:pdfn} \( \mathbf{p}=\mathbf{p}(\vec{v})=
    \tp(\ee{\omega}/M_{\omega})=
    \setof{\varphi(\vec{v})}{\models\varphi[\ee{\omega}]}
    \)\,. 
\end{defn}
    
    Note that since \( \setof{\ee{\alpha}}{\alpha<\olka} \) is
    setwise 
    indiscernible, 
    \( \tupleidx{\ee{\alpha}}{\omega\le\alpha<\olka} \) is a Morley 
    sequence over \( M_{\omega}  \) in \( \mathbf{p} \)\,, the 
    so-called \emph{average type} of \( I_{\omega} \)\,.

\begin{prop}
    \( T \)\nlabel{AlmInd:Tnmd} is non-multidimensional.
\end{prop}

\begin{proof}
    We know that for every ordinal \(\lambda\geq\omega\)\,,
    \(tp(\ee{\lambda}/M_{\lambda})\) is the nonforking extension of
    \(\mathbf{p}\)\,, and moreover \(M_{\lambda}\) is saturated (for
    \(\lambda \geq \tau^{+}\)).  Let \(q\) be over some
    \(M_{\lambda}\)\,, \(\lambda \geq \tau^{+}\)\,.  Set 
    \( \nu=\card{\lambda}^{+} \)\,.  As \(M_{\nu}\) is
    \(\nu\)-saturated, \(q\) is realized in \(M_{\nu}\) which is in
    the algebraic closure of \(M_{\lambda}\) and an independent
    sequence of realizations of \(\mathbf{p}|M_{\lambda}\) (nonforking
    extension of \(\mathbf{p}\) to \(M_{\lambda}\)).  So \(q\) is
    nonorthogonal to \(\mathbf{p}\)\,.  This shows that every type is
    nonorthogonal to \(M_{\omega}\)\,, so \(T\) is
    non-multidimensional.
\end{proof}

So the number of non-orthogonality classes of regular types is 
bounded by \( \tau \)\,, as \( T \) is \( \tau \)-stable and 
\( \card{M_{\omega}}\le\tau \)\,.


\subsection{Structure}\nlabel{AlmInd:Str}

We want to prove a version of \cite[Proposition 2.10]{MR3411166}, and 
explore further consequences of that result. The 
proposition \ref{AlmInd:prp28} generalizing 
\cite[Proposition 2.8]{MR3411166} is essential. In generalizing the 
proofs of \cite[2.8, 2.9, 2.10]{MR3411166} to our more general 
setting, we clarify and improve on many steps of the proofs. For 
background on the ``forking calculus'' arguments, we refer the reader 
to \cite[Chapter 1, \S2]{Pillay:GST}, in particular to 2.20--2.29.

\begin{notation}
Let\nlabel{AlmInd:mub}
\( \mub \) be \( \aleph_{0} \) if \( \mu  \) is finite and 
\( \mu^{+} \) if \( \mu \) is infinite. 
\end{notation}

\tupob Note that \( \mub\le\tau^{+} \)\,.\tupcb

\vsp{1} 

Recall (cf.\ \cite[Lemma 4.4.1]{Pillay:GST}) that \emph{weight 
one} types are the ones for which there is a well-defined dimension 
theory; that non-orthogonality is an equivalence relation on the 
weight one types, and that in particular every regular type has 
weight one (\cite[Lemma 4.5.3]{Pillay:GST}).

Since \( T \) is non-multidimensional (\ref{AlmInd:Tnmd}),
we can make the following definition:

\begin{defn}
   Let\nlabel{AlmInd:OrthCl}
   \( \mathcal{R} \) be the \emph{set} of equivalence classes under 
   non-orthogonality of weight one types of \( T \)\,.  
   
   If \( p \) is some weight one type, \( [p] \) is its class.
\end{defn}

Note that  we
can realize these classes by types  over any \( a \)-model.

\begin{prop} Let\nlabel{AlmInd:prp28}
   \( \lambda\ge\omega \) be an ordinal. 
   
   \tupob There are 
   really only two cases of interest, \( \lambda=\omega \) and
   \( \lambda = \mub \)\,.\tupcb
   
    Consider
    \( M_{\lambda}\prec 
    M_{\lambda+1}=\acl(M_{\lambda}\cup \ee{\lambda}) \)\,.

    There is a set of elements ( 
    \( 
    C=
	 \setof{c_{j}}{j\in J}\subset M_{\lambda+1}\setminus M_{\lambda}
    \)\,, 
	 with \( J \) finite if \( \mu \) is 
    finite and  \( \card{J}\le\mu \) otherwise, such that: 
    \begin{enumerate}
       \item  \( C \) is independent over \( M_{\lambda} \)\,,
    
       \item  each \( \tp(c_{j}/M_{\lambda}) \) is regular,
    
       \item  all regular types occur up to 
       non-orthogonality amongst the types 
       \( \tp(c_{j}/M_{\lambda}) \)\,,
    \end{enumerate}
    and such that \(  M_{\lambda+1} \) is  
    \( a \)-prime and minimal
    over \( M_{\lambda}\cup C \)\,.
    
    Without loss of generality, we can fix some set \( \mathcal{Q} \) 
    of regular types over \( M_{\lambda} \) representing the classes 
    of \( \mathcal{R} \) over \( M_{\lambda} \)\,,
    and assume that for each \( c\in C \)\,,
    \( \tp(c/M_{\lambda})\in\mathcal{Q} \)\,, \tupob so that
    for each \( c,\,c'\in C \)\,,
    either \(\tp(c/M_{\lambda}) = \tp(c'/M_{\lambda})\) or these
   types are orthogonal\tupcb.
\end{prop}

\setcounter{claimno}{0}
\begin{proof}    
    Choose 
    \( 
	 C=\setof{c_{j}}{j\in J}\subset M_{\lambda+1}\setminus M_{\lambda} 
	 \) 
    a maximal independent over \( M_{\lambda} \) set of
    elements realizing regular types over \( M_{\lambda} \)\,.
    Note that by Theorem \ref{AlmInd:Ch}\ref{AlmInd:FaMod}, 
	 \(M_{\lambda} \) is an 
    \( a \)-model. Clearly we can make this choice respecting the 
    final statement of the Proposition.
    \vsp{1}

    \claimnum \(J\) is finite if \(\mu\) is finite and of
    cardinality \(\leq \mu\) otherwise.

   Proof:  This is a weight argument.  In a superstable theory any
   type of a finite tuple \(b\) (over some given base set \(A\)) has
   finite weight in the sense that there is no infinite independent
   over \(A\) set of tuples such that \(b\) forks with each of them
   over \(A\)\,.  For if not, forking calculus gives an
   infinite forking sequence of extensions of \( \tp(b/A)\)\,,
   contradicting superstability.  
   A straightforward extension of this argument shows
   that if \(b\) is a \(\mu\)-tuple then there is no independent over
   \(A\) set of size \(\mu^{+}\) of tuples each of which forks with
   \(b\) over \(A\)\,.  In particular,  each \(c_{j}\)\,, being 
   algebraic over \( M_{\lambda}\cup \ee{\lambda} \) 
   forks with \(\ee{\lambda}\) over \(M_{\lambda}\)\,,
   and so the cardinality
   of \(J\) is at most \( \mu\) when \(\mu\) is infinite.

   \vsp{1}
   \claimnum \(M_{\lambda + 1}\) is
   \( a \)-prime (and therefore prime) and minimal
   over \(M_{\lambda}\cup\setof{c_{j}}{j\in J} \)\,.
   \newline
   Proof.  Let \(N \preceq M_{\lambda + 1}\) be the 
   \( a \)-prime model over
   \(M_{\lambda}\cup\setof{c_{j}}{j\in J} \)\,.  If 
   \(N\neq M_{\lambda + 1}\) then by Lemma \ref{Pre:3mod},
	there is some
   \(d\in M_{\lambda + 1}\setminus N \) whose type
   over \(N\) is regular and does not fork over \( M_{\lambda} \)\,, 
   contradicting the maximality of \(\setof{c_{j}}{j\in J} \)\,.
   
   If \( M_{\lambda + 1} \) is not minimal, 
   then there is \(M_{\lambda}\cup \setof{c_{j}}{j\in J}
   \subseteq N\precneq M_{\lambda+1}\)\,.
   But by Proposition \ref{Pre:ext-a-m} every elementary extension of
   \(M_{\lambda}\) is also an \( a \)-model, so we can repeat the 
   argument just given to get a contradiction.
   
   For the same reason it follows immediately 
   that \(M_{\lambda+1}\) is in fact prime over
   \( M_{\lambda}\cup C \)\,.

   \vsp{1} For the rest, we have already seen (in the proof of 
	\ref{AlmInd:Tnmd}) 
   that any
   regular type \(q\) is nonorthogonal to
   \(tp(\ee{\lambda}/M_{\lambda})\) and so nonorthogonal to a regular
   type
   \(q'\) over \(M_{\lambda}\) which is nonorthogonal to \(p\)\,, and
   so
   realized in the \(a\)-prime model over 
   \(M_{\lambda}\cup \ee{\lambda}\)\,.  But the
   latter is precisely \(M_{\lambda + 1}\)\,.  So \(q'\) is realized
   in
   \(M_{\lambda+1}\)\,, so forks with \(\setof{c_{j}}{j\in J}\) over
   \(M_{\lambda}\)\,.  It easily follows that \(q'\) is nonorthogonal
   to
   some \(tp(c_{i}/M_{\lambda})\)\,.
\end{proof}

\vsp{1} When we say that an infinite tuple \( d \) is \emph{algebraic}
over a set \( A \)\,, we mean that each finite sub-tuple of \( d \) 
is algebraic over \( A \)\,, equivalently that (the range of) \( d \)
is in the algebraic closure of \( A \)\,.

\begin{prop}
    Continuing\nlabel{AlmInd:prp29}
    the notation of Proposition \ref{AlmInd:prp28} 
    \tupop with  \( \lambda=\mub \)\tupcp,
    there are  \( \mu \)-tuples
    \( D=\setof{d_{j}}{j\in J}  \) such 
    that: 
    \begin{enumerate}
        \item  \( \tp(d_{j}/M_{\mub}) \) has weight one and 
        \( c_{j}\in\acl(M_{\mub}\cup\Set{d_{j}}) \) for each 
        \( j\in  J \)\,;
    
        \item \( D  \) is 
        \( M_{\mub} \)-independent; and
    
        \item  \( \ee{\mub} \) is interalgebraic with
        \( D  \) over \( M_{\mub} \)\,.
    \end{enumerate}
    
    \tupob Hence also the types of the \( d_{j} \) represent all the 
    classes of \( \mathcal{R} \) over \( M_{\mub} \)\,.\tupcb
\end{prop}

\setcounter{claimno}{0}
\begin{proof}  
	\( \ee{\mub}=\tupleidx{\ee{\mub,i}}{i<\mu} \)\,. 
   Noting that \( \tp(\ee{\mub}/M_{\mub}) \) is the non-forking 
   extension of \( \mathbf{p} \) to \( M_{\mub} \)\,, for the 
   remainder of this proof we will let \( \mathbf{p} \) denote this 
   type. Set \( C=\setof{c_{j}}{j\in J} \) as given by Proposition
   \ref{AlmInd:prp28}.

   \vsp{1}
   We construct the family \( D \) by a sequence of approximations.
   
   Initially choose \( D \) so that \( D \) is independent over
   \( M_{\mub} \) and for each \( j \)\,, \( d_{j} \) realizes 
   \( \tp(\ee{\mub}/M_{\mub},c_{j}) \)\,.
   
   But then \( D \) is an independent set of realizations of 
   \( \mathbf{p} \) and so \( M'=\acl(M_{\mub}\cup D) \) is a model: 
   an elementary extension of \( M_{\mub} \)\,. In particular, as 
   \( c_{j} \) is algebraic over \( M_{\mub}\cup\Set{\ee{\mub}} \)\,, 
   \( c_{j} \) is algebraic over \( M_{\mub}\cup\Set{d_{j}} \)\,. 
   So \( C \) is contained in \( M' \)\,, and hence 
   \( M_{\mub+1} \) embeds in \( M' \) over 
   \( M_{\mub}\cup C \)\,. Thus  (by taking 
   an automorphism of the universe fixing \( M_{\mub}\cup C \)) we
   can 
   assume without loss of generality that \( M_{\mub+1} \) is 
   contained in \( M' \)\,. Hence:
   
   \claimnum \( \ee{\mub}\in\acl(M_{\mub}\cup D) \)\,.
   
   \vsp{1}
   We now carry out a construction of parameter sequences
   \( f \) which, very informally speaking, encode the domination
   relation between the \( c \)'s and the \( d \)'s. This will 
   eventually allow us to replace each \( d_{j} \) by a 
   \( d_{j}' \) of weight 1 while preserving all the facts proved so 
   far.
   
   Fix \( j \)\,.
   
   \claimnum There is a tuple \(f_{j}\) of length at most \(\mu\) such
   that \(f_{j}\) is independent from \(c_{j}\) over \(M_{\mub}\)
   and \(c_{j}\) dominates \(d_{j}\) over \(M_{\mub}, f_{j}\)\,.
   (That is, if \(a\) is independent from \(c_{j}\) over \(M_{\mub},
   f_{j}\) then \(a\) is independent from \(d_{j}\) over \(M_{\mub},
   f_{j}\)).
   \newline
   Proof.  This is completely standard.  We try to construct a
   sequence \(a_{\alpha}\) of finite tuples, such that such that for
   each \(\alpha\)\,, \(a_{\alpha}\) is independent from \(c_{j}\)
   over
   \(M_{\mub}\cup\setof{a_{\beta}}{\beta < \alpha}\) but
   \(a_{\alpha}\)
   forks with \(d_{j}\) over
   \(M_{\mub}\cup\setof{a_{\beta}}{\beta<\alpha}\)\,.  
   Notice that then for
   each \(\alpha\)\,, \(\setof{a_{\beta}}{\beta \leq \alpha}\)
   is independent
   from \(c_{j}\) over \(M_{\mub}\)\,, but \(d_{j}\) forks with
   \(a_{\alpha}\) over 
   \(M_{\kappa}\cup\setof{a_{\beta}}{\beta < \alpha}\)\,.
   If \(\mu\) is finite then there is (by superstability) a finite
   bound on forking sequences of extensions of
   \(\tp(d_{j}/M_{\kappa})\)\,, and in general, one cannot find such a
   forking sequence of length \(\mu^{+}\)\,.  Hence, for some 
   \(\alpha <\mu^{+}\) one cannot continue the construction to get
   \(a_{\alpha}\)\,.  So take 
   \(f_{j} = \tuple{a_{\beta}: \beta < \alpha}\)\,.

   \vsp{1}
   So for each \( j\in J \) we can choose \( f_{j} \) as described. 
   But by the choice of \( D \) (independent over \( M_{\mub} \)) and 
   the forking calculus we can in fact choose the family of the 
   \( f_{j} \) to be independent from \( C \) over 
   \( M_{\mub} \)\,. Let \( \olc \) denote the concatenation of 
   \( C \) as a \( J \)-tuple, and similarily for 
   \( \old  \) and \( \olf \)\,. Thus we have:
   
   \claimnum \( \olc \) dominates \( \old,\ee{\mub} \) over 
   \( M_{\mub},\olf \)\,, and moreover for each \( j \in  J \)\,, 
   \( c_{j} \) dominates \( d_{j} \) over 
   \( M_{\mub,\olf} \)\,.
   
   \vsp{1}
   We can now find (by superstability and considering the cardinality
   of the relevant set of tuples)  a subset \(A\) of
   \(M_{\mub}\) of cardinality \( \le\tau \) (and therefore 
	\(<\mub \)) such that 
   \(({\olc},{\old}, \ee{\mub}, {\olf})\) is independent from
   \(M_{\mub}\) over
   \(A\)\,.  Let \(A' = \acl^{\Upeq}(A)\)\,, so all types over \(A'\) are
   stationary.

   It follows from the basic facts about forking and domination (see
   \cite[Chapter 1, Lemma 4.3.4]{Pillay:GST}) that:
   \newline
   \claimnum \(\ee{\mub}\in \acl(A', {\old})\)\,,
   \(\tp(c_{j}/A')\) is regular, \(c_{j}\) dominates \(d_{j}\) over
   \(A'\olf\)\,, \(\olc\) dominates \( \old \ee{\mub}\) over
   \(A'{\olf}\)\,, and \({\olc}\) dominates \(\olc \ee{\mub} \) over
   \(A'\) (the latter because \(\olc\) dominates \({\olc}
   \ee{\mub}\) over \(M_{\mub}\)).
   
   \vsp{1} Now using the strong \(\mub\)-saturation of
   \(M_{\mub}\)\,,
   let \(\olf'\in M_{\mub}\) realize \(\tp({\olf}/A')\)\,.  Note
   that \({\olc}\ee{\mub}\) was independent from \({\olf}\) over
   \(A'\)\,, whereby 
   \(\tp({\olf}'{\olc}\ee{\mub}/A') = 
   \tp(\olf\olc \ee{\mub}/A')\)\,.  Now let 
   \({\old}' = (d_{j}':j\in J)\) 
   be such that \(tp({\olf}'{\olc}\ee{\mub},{\old}'/A') =
   \tp({\olf}{\olc}\ee{\mub}{\old}/A')\)\,.  Hence all of Claim 4
   holds with \({\olf}\) replaced by \({\olf}'\)\,, and \(\old\)
   replaced by \({\old}'\)\,:
   
   {\noindent{\emph{Claim} \arabic{claimno}\( {}' \):}\quad}
   \(\ee{\mub}\in \acl(A', {\old'})\)\,,
   \(\tp(c_{i}/A')\) is regular, \(c_{j}\) dominates \(d_{j}'\) over
   \(A'\olf'\)\,, \(\olc\) dominates \( \old \ee{\mub}\) over
   \(A'{\olf'}\)\,, and \({\olc}\) dominates \(\olc \ee{\mub} \) over
   \(A'\)\,.
   
   Note that as \(c_{j}\) dominates \(d_{j}'\) over \(A',{\olf}'\)
   we have that

   \claimnum \(\tp(d_{j}'/A'{\olf}')\) has weight \(1\)\,.

\vsp{1}
Now, as \({\olc}\) is independent from
\(M_{\mub}\) over \(A'{\olf}'\) we have (by domination) that
\({\old}'\) is independent from \(M_{\mub}\) over \(A'{\olf}'\)\,,
and in particular

\claimnum \(\tp(d_{j}'/M_{\mub})\) has weight \(1\) for all
\(i\) (and of course the \(d_{j}'\) are independent over
\(M_{\mub}\)).

\vsp{1} 
Finally, by strong \(\mub\)-saturation of
\(M_{\mub + 1}\) let \({\old}''\) realize 
\(\tp({\old}'/A'{\olf}'{\bar c}\ee{\mub})\) in \(M_{\mub + 1}\)\,.  
Then the domination
statement in Claim 4' implies that \({\old}''\) is independent from
\(M_{\mub}\) over \(A'{\olf}'\) and as in Claim 5, each
\(tp(d_{i}''/M_{\mub})\) has weight 1 and the \(d_{i}''\) are
independent over \(M_{\mub}\)\,.  Moreover \(\ee{\mub}\in acl({\bar
d},M_{\mub})\) (again by Claim 4'), and as \(M_{\mub + 1} =
acl(M_{\mub}, \ee{\mub})\) we conclude that \({\old}''\) is
interalgebraic with \(\ee{\mub}\) over \(M_{\mub}\)\,.  

So we replace the family \( D \) by \( \setof{d_{j}''}{j\in J} \) to 
conclude the proof of the proposition. 
\end{proof}

\begin{thm}
   Any\nlabel{AlmInd:prp210}
   model \(M\) which contains \(M_{\mub}\) is the algebraic
   closure of \(M_{\mub}\) together with an 
   \(M_{\mub}\)-independent set \(D\) of tuples of realizations of
   weight one types over \(M_{\mub}\)\,.
\end{thm}

\begin{proof} 
   Fix some model \(M\) containing \(M_{\mub}\)\,.  Let \( C \) and
   \( \mathcal{Q} \) be chosen over \(M_{\mub}\) as in Proposition 
   \ref{AlmInd:prp28}, for \( \mub \)\,.

   For each \( q\in\mathcal{Q} \)\,,
   let \(I_{q}\) be an enumeration of those \( c\in C \) which realize
   \(q\)\,.  So each \(I_{q}\) is nonempty and
   \(M_{\mub + 1}\) is \( a \)-prime and minimal over
   \(M_{\mub}\cup\bigcup_{q}I_{q}\)\,.  
   Let \(\lambda_{q}\le\mub \) be the cardinality of
   \(I_{q}\)\,.
   
    For some sufficiently large
   \(\lambda\)\,, let \(M'\) be a saturated model of cardinality
   \(\lambda\) extending \(M\)\,.

   Since \(M\) extends \( M_{\mub} \)\,, \( M \) is an 
   \( a \)-model.   For each \( q\in\mathcal{Q} \) let 
   \(I_{q}'\) be a
   maximal \(M_{\mub}\)-independent set of realizations of \(q\)
   in \(M\)\,. Note that it may be empty.
   Let \( C'=\bigcup_{q\in\mathcal{Q}}I_{q}' \)\,.
   As in the proof of Proposition \ref{AlmInd:prp28}, \(M\) is
   \( a \)-prime and minimal over 
   \(M_{\mub}\cup C' \)\,.  Now we can extend
   each \(I_{q}'\) to \(I_{q}''\)\,, a maximal
   \(M_{\mub}\)-independent
   set of realizations of \(q\) in \(M'\)\,.  Then \(I_{q}''\) has
   cardinality \(\lambda\)\,. Let 
   \( C''=\bigcup_{q\in\mathcal{Q}}I_{q}'' \)\,.

   For each \(c\in \bigcup_{q}I_{q}''\) we can 
   (by Proposition \ref{AlmInd:prp29})
   choose suitable \(d_{c}\) whose type over \(M_{\mub}\) is of
   weight \(1\) and with \(c\) dominating \(d_{c}\) over
   \(M_{\mub}\) and  \(c\in \acl(M_{\mub},d_{c})\)\,, such that
   \begin{enumerate}
      \item if \(c\in M\) (i.e. \(c\in I_{q}' \) for some 
      \(q\in\mathcal{Q}\)) then
      \(d_{c}\in M\)\,, and
   
      \item  \(M' = \acl(M_{\mub} \cup \bigcup_{c\in C''}d_{c})\)\,.
      
      \vsp{1}
      Let \(B = \acl(M_{\mub} \cup\setof{d_{c}}{c\in C'}\)\,.  
      So:
   
      \item  
      \(
      M_{\mub}\cup C' \subseteq B \subseteq M
      \)\,.
   \end{enumerate}

   As  \( \setof{d_{c}}{c\in C''} \) is independent
   over \(M_{\mub}\)\,, it follows that for any \(c_{1},..,c_{r}\in
   C''\setminus C'\)\,,
   \(\tp(d_{c_{1}},...,d_{c_{r}}/B)\) is finitely satisfiable in
   \(M_{\mub}\)\,.  It follows that \(B\) is the universe of an
   elementary substructure of \(M'\)\,, so by (3) and the minimality
   of \(M\) over \(M_{\mub}\cup C'\) we have that \(M = B\)\,,
   proving the proposition.
\end{proof}

\begin{cor}
   Continuing\nlabel{AlmInd:hulls}
   the notation of the preceding results, we can find
   \tupop in \( M_{\mub+1} \)\tupcp\ sets 
   \( \setof{N(r)}{r\in\mathcal{R}} \)\,, each
   uniquely determined up to isomorphism over \( M_{\mub} \) by
   \( r \)\,, with \( \tp(N_{r}/M_{\mub})\in r \) and
   such that each \( N(r) \) is  a maximal \tupop with respect to 
	\( \subseteq \)\tupcp\ weight one set over
   \( M_{\mub} \)\,.
   
   We call \( N(r) \) the \emph{hull} of \( r \) 
   \tupop over \( M_{\mub}  \)\tupcp. 
   
   Furthermore, if \( M \) is any model containing \( M_{\mub} \)
   \tupop as in \ref{AlmInd:prp210}\tupcp, then \( M \) is the 
   algebraic closure of a
   family \tupop independent over \( M_{\mub} \)\tupcp\ 
   of copies of the various \( N(r) \)\,, 
   \( r\in \mathcal{R}  \)\,.
\end{cor}
 
\begin{proof}
   For each \( r\in R \)\,, let \( c \) realize a type in
   \( \mathcal{R} \) over \( M_{\mub} \)\,, and let \( N(r) \) be the
   \( a \)-prime model over \( M_{\mub}\cup\Set{c} \)\,. The type  
	\( q=\tp(N(r)/ M_{\mub} ) \) is dominated by 
	\( \tp(c/M_{\mub}) \)\,, so \( q \) is weight one and in the class 
	\( r \)\,. If \( X \) is any set whose type over  \( M_{\mub} \) is 
	of weight one and in \( r \) then the \( a \)-prime model over
	\( M_{\mub}\cup X \) is isomorphic to \( N(r) \)\,, and so 
	\( N(r) \) is maximal with the stated properties, and unique up to 
	isomorphism.
	
    Note that \( N(r)\supset M_{\mub} \)\,.

   Furthermore, if \( M \) is any model containing \( M_{\mub} \)
   \tupop as in \ref{AlmInd:prp210}\tupcp, then as in the proof of
   that theorem we find
   \( M= \acl(M_{\mub} \cup\setof{d_{c}}{c\in C'}) \) where
   \( \setof{d_{c}}{c\in C'} \) is an \( M_{\mub} \)-independent set
   of tuples realizing weight one types over \( M_{\mub} \)\,.
   Let \( r_{c}=[\tp(d_{c}/M_{\mub})] \)\,.
   Then for each \( c\in C' \)\,, \( d_{c} \) is contained in a copy
   \( N_{c} \)
   of \( N(r_{c}) \)\,, which by definition is domination-equivalent 
	to \( d_{c} \)\,, and so we can choose
   \( \setof{N_{c}}{c\in C'} \) to be independent over
   \( M_{\mub} \)\,. Then of course
   \( M=\acl\left(\bigcup_{c\in C'}N_{c}\right) \)\,.
\end{proof}

   In particular cases, the structure theory can be refined quite a 
   bit.

\begin{exm} The primordial motivating example for stability theory 
	 already exhibits this structure, and more.
   Let \( T \) be the theory of algebraically closed fields of some 
   fixed characteristic. Then a transcendence basis is an 
   indiscernible set, and any model is the algebraic closure of its 
   transcendence basis. So \( T \) is \( (1,\aleph_{1}) \)-almost 
   indiscernible.
   
   The results of this section describe a structure theory for the 
   extensions of the model with countably many independent 
   transcendental elements (\( \mub = \aleph_{0} \)), but of course
    here we 
   actually have a structure theory for extensions of the prime model.
   
   \( \mathcal{R} \) consists of a unique class, and the hull of that 
   class is the field of transcendence degree one. So 
	Corollary \ref{AlmInd:hulls} sees every 
   algebraically closed field as the algebraic closure of an 
   algebraically independent family of algebraically closed fields of 
   trancendance degree one.
\end{exm}

   In the next subsection, we will see similar kinds of examples in 
   theories of modules.
   

 \subsection{The case of theories of modules}\nlabel{ssMod:}

   We take Prest's book \cite{Pr:book} as our main reference, to 
   ensure a uniform approach to the subject, with occasional 
   attributions to  primary sources.
   
   Throughout, \( T \) is a \textbf{complete superstable theory}
   of \( R \)-modules, \( \card{T}=\tau \)\,.  

   \( \lambda(T) \) is the least cardinal in which \( T \) is stable
   (so \( \tau\le\lambda(T)\le 2^{\tau} \)).

\vsp{1}
\begin{prop}
    \upshape{(\cite[Cor 3.8]{Pr:book}, due to Ziegler)}
   If\nlabel{ssMod:tt} \( M\prec N\models T \) then the factor module 
	\( N/M \) is totally transcendental.
\end{prop}

We learn the following facts from Prest \cite[\S 6.5]{Pr:book} 
(originally Pillay-Prest \cite{PiPr:ModSt}):

\begin{prop}{\quad}
    
   \begin{enumerate}
      
      \item \( M\models T \)\nlabel{ssMod:aMod}
      is an \( a \)-model if{}f \( M \) is
      pure-injective and weakly saturated. 
		{\upshape\cite[6.37]{Pr:book}}

      \item  Pure-injective models of \( T \) are
      \( \aleph_{0} \)-homogeneous. 
		{\upshape\cite[6.35]{Pr:book}}

      \item Elementary extensions of pure-injective models are 
      pure-injective. 
		{\upshape\cite[6.34]{Pr:book}}

      \item Elementary extensions of \( a \)-models are 
      \( a \)-models. 
		{\upshape\cite[6.41]{Pr:book}}
   
   \end{enumerate}
\end{prop}

\begin{prop}
   In\nlabel{ssMod:api}
   general, for any complete theory of modules \( T' \)\,, if every 
   \( a \)-model is pure-injective, then \( T' \) is superstable.
	{\upshape\cite[6.40]{Pr:book}}
\end{prop}

\begin{exm}
   Consider\nlabel{ssMod:p-adic}
   the \( p \)-adics  \( \overline{\ZZ_{(p)}} \)  in two
   ways, as a \( \ZZ \)-module and as 
	\( \overline{\ZZ_{(p)}} \)-module.  In both
   cases the theory is superstable not totally transcendental, with 
   \( \lambda(T)=2^{\aleph_{0}} \)\,.  But the latter has 
   \( \lambda(T)=\card{T} \)\,. (These are used as illustrative 
	examples of many aspects of the model theory of modules throughout 
	\cite{Pr:book} and the facts stated here are ``common knowledge''.)
   
   Models of the theory as a \( \ZZ \)-module have the form 
   \( M\oplus \QQ^{(\kappa)} \)\,, where 
   \( M\preceq \overline{\ZZ_{(p)}} \) and \( \kappa\ge 0 \) is a 
   cardinal.  Models of the theory as a 
   \( \overline{\ZZ_{(p)}} \)-module have the form 
   \( \overline{\ZZ_{(p)}}\oplus {\QQ_{p}}^{(\kappa)}\)\,, where 
   \( {\QQ_{p}} \) is the quotient field.
   
   In either case, there are no algebraic or definable elements other 
   than \( 0 \)\,. Note however that the type of, for instance,
   \( 1\in \overline{\ZZ_{(p)}} \)\,, while not algebraic, is 
   \emph{limited} in the sense that the pure-injective hull of a 
   realization of it occurs exactly once as a direct summand of any 
   model of \( T \)\,. As a  \( \ZZ \)-module, the theory is not 
   \( (\mu,\kappa) \)-almost indiscernible for any countable 
   \( \mu \)\,. But as a \( \overline{\ZZ_{(p)}} \)-module, the
   \( 2^{\aleph_{0}} \)-saturated model 
   \( \overline{\ZZ_{(p)}}\oplus {\QQ_{p}}^{(2^{\aleph_{0}})}\) is the 
   \emph{definable} closure of an indiscernible set of tuples of 
   cardinality \( 2^{\aleph_{0}} \)\,, where for convenience we take 
   the order type of the tuples to be 
   \( 2^{\aleph_{0}} + 2^{\aleph_{0}} \)\,, the first 
   \( 2^{\aleph_{0}} \) components of the tuple being some fixed 
   enumeration of \( \overline{\ZZ_{(p)}} \)\,, and the second 
   sequence of \( 2^{\aleph_{0}}  \) components ranging over an 
   enumeration of the standard basis for 
   \( {\QQ_{p}}^{(2^{\aleph_{0}})} \)\,. (The details of this 
   construction are 
   made explicit in the proof of Theorem \ref{AlmInd:ssMod} 
   following.)
\end{exm}
   
   In fact, the situation described in this example is typical:

\begin{thm}\nlabel{AlmInd:ssMod}
   Let \( T \) be a superstable theory of modules with 
   \( \lambda(T)=\card{T}=\tau \)\,, \tupop in particular,  if 
   \( T \) is totally transcendental\tupcp.  Then \( T \) is an almost
   indiscernible theory of modules.
\end{thm}

\begin{proof}
   We extract the required properties of \( a \)-models in 
   superstable theories of arbitrary cardinality from 
   Baldwin's book \cite{Bald:FunStb} on stability theory, 
    Chapter XI, \S1, \S2. 
   
   In particular, there is an \( a \)-prime model
   \( \mathcal{A}_{0} \) 
   (of cardinality \( \tau \)). There is a saturated 
   proper elementary extension 
   \( N \) of \( \mathcal{A}_{0} \) of cardinality  
   \( \tau^{+} \)\,, 
   since \( T \) is stable in all cardinals greater than
   \( \tau \)\,.  
   Any type 
   that is realized in \( N\setminus \mathcal{A}_{0} \) 
   is realized by an 
   independent set of cardinality \( \tau^{+} \)\,.
   The factor module
	\( N/\mathcal{A}_{0} \) is a tt module by \ref{ssMod:tt}, which 
   decomposes as direct sum of indecomposables. In particular, since 
   \( N \) and \( \mathcal{A}_{0} \) are themselves pure-injective,
   each 
   summand is a direct summand of \( N \)\,. These summands are 
   necesarily of cardinality (less than or equal to) \( \tau \)\,, 
   the 
   cardinality of the language, as the theory of
   \( N/\mathcal{A}_{0} \) is tt, and there are, 
   up to isomorphism, no more than \( \tau \) distinct summands. 
   There are no
    ``limited'' summands (summands which appear a fixed finite number 
    of times in any model) as these all necessarily appear as 
    summands of \( \mathcal{A}_{0} \)\,. Therefore (since \( N \) has
    a large 
    independent set over \( \mathcal{A}_{0} \)) every summand of 
    \( N/\mathcal{A}_{0} \) 
    occurs \( \tau^{+} \) times. Let \( A \) be the direct sum of one 
    copy, up to isomorphism, of each summand of 
	 \( N/\mathcal{A}_{0} \)\,. By 
    the arguments just given, \( \card{A}\le\tau \)\,.
    Then \( N\cong \mathcal{A}_{0}\oplus A^{(\tau^{+})} \)\,. 
    Just as we did in Example \ref{ssMod:p-adic}, fix an enumeration 
    \( \overline{m} \) of \( \mathcal{A}_{0} \) in order type 
	 \( \tau \) and 
    an enumeration \( \ola \) of \( A \) of order type 
    \( \le\tau \)\,; for \( i<\tau^{+} \) let \( \ola_{i} \) be the 
    copy of \( \ola \) on the \( i \)-th component of
    \( A^{(\tau^{+})} \)\,. Let \( \ee{i}  \) be the concatenation of 
    \( \overline{m} \) and \( \ola_{i} \)\,.
    Clearly \( \setof{\ee{i}}{i<\tau^{+}} \) is a set of sequences
    all 
    of the same type and independent, since direct-sum independent, 
    and so is a set of indiscernibles.
    
    Thus \( T \) is seen to be \( (\tau,\tau^{+}) \) almost 
    indiscernible.
\end{proof}

\begin{remnc}
   In the case where \( T \) is tt, we can carry out the construction
   just described, taking \( \mathcal{A}_{0} \) to be the sum of the
   limited summands of \( T \)\,, if there are any, or \( 0 \) 
   otherwise. So the choice of \( \mathcal{A}_{0} \) as described in 
   the proof of the Theorem does not necessarily give the sharpest 
   possible structure theorem. Nor will this crude construction 
   reveal whether or not \( T \) is \( (\mu,\tau^{+}) \) almost 
   indiscernible for some \( \mu<\tau \)\,.
\end{remnc}

\begin{cor}
    A\nlabel{AlmInd:AIss}
    complete theory \( T \) of modules is almost indiscernible 
    if{}f it is superstable with  \( \lambda(T)=\card{T} \)\,.
\end{cor}

\begin{proof}
    By Theorems \ref{AlmInd:Tss} and \ref{AlmInd:ssMod}.
\end{proof}

\begin{exm}
    Consider\nlabel{AlmInd:RRtt}
    the following example, used at several places in 
    Prest \cite{Pr:book}, in particular at 2.1/6(vii)
    (with \( k \)  a countable field). This was an 
    important example of Zimmermann-Huisgen and Zimmermann
    \cite{MR0498722}\,.
    
    Let \( k \) be an infinite field (of cardinality \( \tau \)) and
    set
    \( 
    R=k[(x_{i})_{i\in\omega}:x_{i}x_{j}=0\mbox{ for all }i,j]
    \)\,.  
    (We can of course, with some small adjustments to the
    cardinalities in what follows, make the same construction with an
    uncountable family of indeterminates.)  Then \( \lsub{R}{R} \) is
    an indecomposable tt module.  Its lattice of pp definable
    subgroups consists of all the finite dimensional vector subspaces
    of \( J=\bigoplus_{i\in\omega}kx_{i} \)\,, together with \( J \)
    itself and \( R \)\,.  Morley rank equals Lascar rank.  Each
    subspace of dimension \( n \) has rank \( n \)\,; \( J \) has rank
    \( \omega \)\; and \( R \) has rank \( \omega+1 \)\,.  There are
    thus two indecomposable pure-injective summands of models of 
    \( T \)\,: \( \lsub{R}{R} \) corresponding to the types of finite
    rank/rank \( \omega+1 \)\,, and \( \lsub{R}{k} \) corresponding to
    the type of rank \( \omega \)\,.  In the latter case, the action
    of \( R \) on \( k \) is given by \( x_{i}k=0 \) for all 
    \( i \)\,.
    
    The models of the theory of \( \lsub{R}{R} \) are precisely the 
    modules \( {\lsub{R}{R}}^{(\kappa)}\oplus k^{(\lambda)} \)\,, 
    with \( \kappa\ge 1 \) and \( \lambda\ge 0 \)\,. For 
    \( \kappa\ge\tau \)\,, the 
    saturated model of power \( \kappa \) is 
    \( {\lsub{R}{R}}^{(\kappa)}\oplus k^{(\kappa)} \)\,.
    
    So in particular although each free module 
    \( {\lsub{R}{R}}^{(\kappa)} \)\,, \( \kappa \) a (non-zero) 
    finite or infinite cardinal, is in the algebraic closure of 
    an indiscernible set of cardinality \( \kappa \) (just take the 
    standard basis vectors), these models are not saturated.
    
    However of course this theory is \( (2,\tau^{+}) \)-almost 
    indiscernible, since each indecomposable is 1-generated.
\end{exm}


\section{Free algebras}\nlabel{Free:}


\subsection{The general theory}\nlabel{Free:Th}

We provide generalizations and extensions of the results in Pillay and
Sklinos \cite[Section 3]{MR3411166} to the uncountable context.  None
of the proofs depended in any significant way on the assumption that
the language is countable, and go through with \( \aleph_{1} \)
replaced by \( \tau^{+} \)\,.  But we verify all the details in any
case. As in \cite{MR3411166}, we refer the reader to the text
 \cite{MR648287}, \cite{BurSan} 
of Burris and Sankappanavar for the elementary facts of 
Universal Algebra.

If \( N \) is an algebra
and \( X\subseteq N \)\,, then \( \gen{X} \) is the subalgebra
of \( N \) generated by \( X \)\,.

We start off with a couple of simple facts about free algebras in a 
variety. Note that the cardinality of a free basis is \emph{not} in 
general an 
invariant of a free algebra unless that cardinality is greater than 
the cardinality of the language, c.f.\ Example \ref{Free:Mod}.

\begin{lemma}
    {\upshape\cite[Remark 3.1]{MR3411166}}\nlabel{Free:bas}
    Suppose the algebra \( A  \) is free on \( X_{1}\cup X_{2} \)\,,
    \( X_{1} \cap X_{2}=\emptyset \)\,.
    Let \( A_{1}=\gen{X_{1}} \)\,, so \( A_{1} \) is free on
    \( X_{1} \)\,. Let \( Y_{1} \) be any other free
    basis for \( A_{1} \)\,. Then \( A \) is free on 
    \( Y_{1}\cup X_{2} \)\,.
\end{lemma}

We give here a more general version of \cite[Lemma 3.7]{MR3411166},
which is a fact of universal algebra, not a consequence of the 
context in which we work. 

\begin{lemma}
    Let{\upshape\nlabel{Free:37}} \( M \) be a free algebra on
    \( \tau^{+} \) generators in a 
    variety \( \mathcal{V} \) over a language of cardinality 
    \( \tau \)\,. Suppose that \( Y \) is a free basis of \( M \)\,; 
    \( \ola \) is a finite tuple in some other free basis \( X \) of 
    \( M \)\,; and for some finite tuple \( \oly\in Y \) and finite 
    tuple of terms \( \overline{t} \)\,, 
    \( \ola=\overline{t}(\oly) \)\,.

    Then for any 
    \( C\subset Y\setminus\oly \) such that 
    \( \card{Y\setminus C}=\tau^{+} \)\,,
    \( C\cup\Set{\ola} \) may be extended to a basis of \( M \)\,.
\end{lemma}

\begin{proof}
    Let \( Y_{0}=Y\setminus C \)\,.  
    So \( \card{Y_{0}}=\tau^{+} \)\,, \( \oly\in Y_{0} \)\,,  and 
    \( \ola\in\gen{Y_{0}} \)\,.
    Then \( \gen{Y_{0}}\cong M \) as they 
    are both free on \( \tau^{+} \) generators; 
    \( \gen{Y_{0}}\subset M \) (and is a proper subset if \( C \) is 
	 nonempty);  and
    \( C\cap \gen{Y_{0}}=\emptyset \) as \( Y \) is a free basis.
    So by Lemma \ref{Free:bas}, \( C\cup\Set{\ola} \) may be extended 
    to  a basis of \( M \)\,.
\end{proof}

So in particular if \( a \) belongs to a basis \( X \) of \( M \)\,,
and \( Y  \) is a 
basis of \( M \)\,, then there is \( b\in Y \) such that 
\( \Set{a,\,b} \) can be extended to a basis of \( M \)\,.

\vsp{1}
\begin{context}
    Take\nlabel{Free:Ctx}
    as a replacement for 
    {\upshape\cite[Assumption 3.2]{MR3411166}} the 
    following:
    
    Let \( \mathcal{V} \) be a variety over an algebraic language
    \( \LL \) of 
    cardinality \( \tau\ge\aleph_{0} \)\,.
    Let the algebra \( M \) be a free algebra for 
    \( \mathcal{V} \) on a 
    set \( I=\setof{\ee{\alpha}}{\alpha<\tau^{+}} \) 
    \tupop of individual 
    elements\tupcp, such that \( M \) is \( \tau^{+} \) saturated.
\end{context}

Adopt the same notational conventions as in Context \ref{AlmInd:Ctx}. 
So in particular the underlying theory \( T \) 
is the theory of \( M \)\,. 

\begin{lemma}\quad\nlabel{Free:extn}
    \begin{enumerate}
        \item  \( I \) is a set of indiscernibles in 
        \( M \)\,.
    
        \item  If \( I'\subset I \)\,, or if
        \( I'\supset I \) and \( I' \) is a set of indiscernibles in 
        \( \Mons \) extending \( I \)\,, then
        \( \gen{I'} \) is  free on \( I' \) in \( \mathcal{V} \)\,.
    \end{enumerate}
\end{lemma}

\begin{proof}
    (1) holds in general by freeness; 
    the first case of (2) always holds for 
    a free algebra, and the second case of (2) then follows by 
    indiscernibilty and
    the homogeneity of the universe. For clearly then any subset 
    \( I_{0}\subset I' \)
    of cardinality \( \tau^{+} \) is a free basis for
    \( \gen{I_{0}} \)\,; the set of all such subsets of 
    \( I' \) is an updirected family; if 
    \( f:I'\rightarrow A \) for some algebra
    \( A\in\mathcal{V} \) then each \( f\restr I_{0} \) has a 
    unique lifting to a homorphism 
    \( \overline{f_{I_{0}}}\rightarrow A \)\,, and these liftings are 
    pairwise compatible (else we would have a contradiction to 
    indiscernibility). So the union of the maps
    \( \overline{f_{I_{0}}} \) lifts \( f \) to \( A \)\,.
\end{proof}

\begin{cor}\quad\nlabel{Free:ss}
    
    \begin{enumerate}
        \item  \( T \) is superstable.
    
        \item  If \( T \) is a theory of modules, then 
        \( T \) is tt, and \( T=\Th(F^{(\aleph_{0})}) \)\,, where 
        \( F \) is the free module in \( \mathcal{V}  \) on one 
        generator.
    \end{enumerate}
\end{cor}

\begin{proof}\quad
    
    \begin{enumerate}
        \item  By Theorem \ref{AlmInd:Tss}.
    
        \item  The free module on a set \( I \) in 
        \( \mathcal{V} \) is isomorphic to 
        \( F^{(I)} \)\,, where \( F \) is the free module in 
        \( \mathcal{V} \) on one generator. So \( T \) is a 
        superstable theory of modules closed under products, hence is 
        tt. All infinite weak direct powers of any module are 
        elementarily equivalent, 
        so \( F^{(I)}\equiv  F^{(\aleph_{0})}  \)\,.
    \end{enumerate}
\end{proof}

\begin{defn}
    Following\nlabel{Free:bsc}
    {\upshape\cite[Definition 3.3]{MR3411166}}, we call 
    \( B\subset M \) \emph{basic} if it is a subset of some free 
    basis of \( M \)\,. We call \( b\in M \) \emph{basic} if
    \( \Set{b} \) is basic.
\end{defn}

\begin{lemma}{\upshape\cite[Lemma 3.4]{MR3411166}}
    There is\nlabel{Free:p0}     
    a type \( p_{0} \) over \( \emptyset \) such that for any
    \( a\in M \)\,, \( a \) is basic if{}f \( a \) realizes 
    \( p_{0} \)\,.
\end{lemma}

\begin{proof}
    Since \( I \) is an indiscernible set, all the elements of 
    \( I \) have the same type \( p_{0} \)\,. If \( X \) is any other 
    basis of \( M \)\,, since \( \card{M}>\tau \)\,, 
    \( X \) also has cardinality \( \tau^{+} \)\,, and any bijection 
    between \( X  \) and \( I \) extends, by freeness, to an 
    automorphism of \( M \)\,. So the elements of \( X \) also have 
    type \( p_{0} \)\,.
    
    Conversely, if the type of \( a \) is \( p_{0} \) and 
    \( e\in I \) then by saturation there is an automorphism
    \( f \) of \( M \) taking \( e \) to \( a \)\,; then 
    \( f[I] \) is a free basis containing \( a \)\,.
\end{proof}

\vsp{1}
\textbf{Fix} \( p_{0} \) as in Lemma \ref{Free:p0}.

\begin{remnc}
    We would like to see that the rank of \( p_{0} \) is maximal. All 
    that is needed is that if \( M \)\,, \( N \) are models of a 
    superstable theory, and \( f:M\rightarrow N \) is a homomorphism, 
    then  the \( U \)-rank (Morley rank, as the case may be) of 
    \( a\in M \) is greater or equal the \( U \)-rank (Morley rank) 
    of \( f(a) \)\,. At present, we only have the result for theories 
    of modules.
\end{remnc}

\begin{lemma}
    {\upshape\cite[Lemma 3.8]{MR3411166}}\nlabel{Free:38}
    The type \( p_{0} \) is stationary. Hence so is
    \( {p_{0}}^{(n)} \) for any \( n \)\,.
\end{lemma}

\begin{proof}
    We have to show that \( p_{0} \) determines a unique strong type 
    over \( \emptyset \)\,. 
    
    Suppose that \( a \) and \( b \) are realizations of 
    \( p_{0} \)\,. So \( a \) is an element of some basis 
    \( X \) of \( M \) and \( b \) is an element of some basis
    \( Y \) of \( M \)\,. By Lemma \ref{Free:37}, there is
    \( b'\in Y \) such that \( \Set{a,\,b'} \) is basic, so extends to 
    a basis \( Z \) of \( M \)\,.  But \( Z \) is indiscernible in 
    \( M \)\,, so \( a  \) and \( b' \) have the same strong type.
    But \( b \) and \( b' \) have the same strong type (as elements 
    of \( Y \)), so \( a  \) and \( b \) have the same strong type.
\end{proof}

\begin{cor}
    All\nlabel{Free:Stat} \( n \)-types over \( \emptyset \) are 
    stationary.
\end{cor}

\begin{proof}
    [From the proof of \cite[Corollary 3.10, 3.11]{MR3411166}]
    If \( a\in M \) then \( a\in\dcl(e) \) for some finite sequence
    \( e\in I \)\,. But \( \tp(e/\emptyset) \) is stationary by 
    Lemma \ref{Free:38}. 
    Therefore \( \tp(a/\emptyset) \) is stationary.
\end{proof}

\begin{cor}
    {\upshape\cite[Corollary 3.10]{MR3411166}}\nlabel{Free:310}
    \( \acl^{\Upeq}(\emptyset)=\dcl^{\Upeq}(\emptyset) \)\,.
\end{cor}

\begin{proof}
    Immediate.
\end{proof}

\begin{prop}
    {\upshape\cite[Proposition 3.9]{MR3411166}}\nlabel{Free:39}
    The sequence \( I \) is a Morley sequence in \( p_{0} \)\,.
\end{prop}

\begin{proof}
    We have to show that \( I \) is independent over 
    \( \emptyset \)\,. 
    
    Let \( I_{0}=\setof{\ee{n}}{n<\omega} \)\,. By homogeneity it is 
    enough to  show that \( \ee{\omega} \) and \( I_{0} \) are 
    independent over \( \emptyset \)\,. Let \( a \) realize a 
    non-forking extension of \( p_{0} \) to \( I_{0} \)\,.
    By Lemma \ref{Free:37} there is an infinite 
    \( I_{0}'\subseteq I_{0} \) such that 
    \( I_{0}'\cup{a} \) is basic. But then by freeness there is an 
    automorphism carrying \( I_{0}'\cup{a} \) to 
    \( I_{0}'\cup{\ee{\omega}} \)\,, so \( \ee{\omega} \) and 
    \( I_{0} \) are independent.
\end{proof}

\begin{prop}
    {\upshape\cite[Proposition 3.12]{MR3411166}}\nlabel{Free:312}
    Let \( \ola,\,\olb\in M \)\,.
    
    Then \( \ola \) is independent from \( \olb \) over 
    \( \emptyset \) if{}f there is a basis \( A\cup B \)\,,
    \( A,\, B \) disjoint, of
    \( M \) such that \( \ola\in \gen{A} \) and 
    \( \olb\in \gen{B} \)\,.
\end{prop}

\begin{proof}
    The reverse direction is clear by Proposition \ref{Free:39}, as
    any 
    basis is an independent set.
    
    For the forward direction, suppose 
    \( \ola \) is independent from \( \olb \) over 
    \( \emptyset \)\,. Without loss of generality, 
    for some \( n<\omega \)\,, 
    \( \ola, \olb\in\gen{\ee{i}:i<n} \)\,. In particular,  
    \( \ola \) is expressed as a \emph{sequence} of terms 
    in \( \Set{\ee{i}:i<n} \)\,, 
    \( \ola=\vec{t}(\ee{i}:i<n) \)\,.
    Let \( \ola'=\vec{t}(\ee{i}:n\le i<2n) \)\,. Then
    \( \ola' \) is independent from \( \olb \) (by the reverse 
    direction already proved!) and 
    \( \tp(\ola)=\tp(\ola') \)\,.
	 Thus there is an 
    automorphism of the universe fixing \( \olb \) and taking 
    \( \ola' \) to \( \ola \)\,. The image of the basis \( I \) under 
    this automorphism gives us the required decomposition.
\end{proof}

The proof extends in the obvious way to infinite tuples, and to 
independence over an arbitrary basic set.


\subsection{The particular case of modules}

 Let
 \( \mathcal{V} \) be a variety of (left) \( R \)-modules.  The
 free module \( \mathcal{F}_{1} \) on one generator in 
 \( \mathcal{V} \) is clearly an image of the 
 (absolutely) free module on one
 generator, that is, of \( \lsub{R}{R} \)\,.  We take \( 1\in R \)
 as the free generator.  So \( \mathcal{F}_{1}\cong R/I \) for some
 left ideal \( I \) of \( R \)\,.  The free module in 
 \( \mathcal{V} \) on a set \( X \) is then (up to isomorphism)
 \( \mathcal{F}_{1}^{(X)} \)\,.
 In particular, by ``the free module in \( \mathcal{V} \) 
  on \( \kappa \)-many
 generators'', we mean \( \mathcal{F}_{1}^{(\kappa)} \)\,.

\begin{exm}
    Note\nlabel{Free:Mod}
    however that free modules on different cardinalities need not
    be distinct: take \( K \) a field, \( \kappa \) an infinite 
    cardinal,
    and set \( R \) to be the ring
    of all column-finite \( \kappa\times\kappa \) ``matrices'' over 
    \( K \)\,.  It is an easy exercise to verify that the usual matrix
    multiplication is well-defined.  By considering partitions of
    \( \kappa \) into \( 2,\,3,\,\ldots,n,\ldots \) pairwise disjoint
    subsets, each of cardinality \( \kappa \)\,, we see that
    \[
    \lsub{R}{R}\cong{\lsub{R}{R}}^{(2)}\cong
    {\lsub{R}{R}}^{(3)}\cong\,\cdots\,\cong 
    {\lsub{R}{R}}^{(n)}\cong\cdots 
    \]
    However the cardinality of a free basis is uniquely defined 
	 whenever it is infinite.
\end{exm}

 Let \( \kappa\ge\card{R}^{+} \)\,.
 Suppose that \( M=\mathcal{F}_{1}^{(\kappa)} \)\,, the free module
on 
 \( \kappa \) generators in \( \mathcal{V} \)\,, is saturated and 
 let \( T=\Th(M) (=\Th(\mathcal{F}_{1}^{(\aleph_{0})})  \)\,. Since 
 \( T  \) is a superstable theory of modules and closed under 
 products, it is tt. 

Recall the ordering of pp-types of a theory \( T \) of modules: 
if \( p \) and \( q \) are pp-types of \( T \)\,, then
\( p\le q \) if for all \( M\models T \)\,,
\( p[M]\subseteq q[M] \)\,. Equivalently, \( p\le q \) if{}f
\( p\supseteq q \)\,. When we say ``maximal'', we mean 
``maximal with respect to the ordering \( \le \) on pp-types''.

\begin{prop}
    The\nlabel{Free:Mod:max} pp-type \( {p_{0}}^{+} \) 
    of a basic element
    is maximal,  hence  
    \( p_{0} \) has maximal Morley rank. 
\end{prop}

\begin{proof}
    Fix a basis \( X \) of \( M \) and \( e\in X \)\,. 
    So the pp-type of \( e \) is \( {p_{0}}^{+} \)\,.
    \( N\equiv M \) certainly implies \( N\in\mathcal{V} \)\,.
    Let \( q \) be a pp-type of \( T \) and let
    \( a\in N\models T \) realize \( q \)\,. Define
    \( f:X\rightarrow N \) by setting \( f(e)=a \) and letting
    \( f \) be arbitrary otherwise. Then by freeness 
    \( f \) extends to a homomorphism 
    \( \overline{f}:M\rightarrow N \)\,, and homomorphisms increase 
    pp types setwise, so \( {p_{0}}^{+}\subseteq q \)\,, that is,
    \( {p_{0}}^{+}\ge q \)\,.
\end{proof}

\vsp{1}
Pillay and Sklinos ask (\cite[Question 3.14]{MR3411166}) 
whether the theory of a 
saturated free algebra must have finite Morley rank, and suggest that 
the answer should be easy in a variety of \( R \)-modules. 
Under suitable restrictions, the answer is indeed ``yes''. For 
instance,
if \( \mathcal{V} \) is a variety of \( R \)-modules
such that the free module \( \lsub{R}{N_{1}} \) on one generator has a
unique indecomposable direct summand, and the free module 
\( \lsub{R}{N} \) on \(\card{R}^{+} \) generators is saturated, then 
\( \Th(N) \) is unidimensional and so has finite Morley rank. 
However the answer is ``no'' in general; see Example 
\ref{Free:ModUTMR} following.

On the other hand it is natural to ask if there is classification of
those rings \(R\) such that free \(R\)-modules on \(|R|^{+}\) 
generators are saturated. (In fact this was explicitly asked by Piotr
Kowalski during a talk on the subject by the second 
author in Wroclaw, Poland, in July 2019.) This is answered in the
theorem below.  See Chapter 14 of \cite{Pr:book} for 
definitions of the notions of coherence and perfectness
of a ring \(R\)\,.  The theorem and proof are thematically 
close to the material in this Chapter; a 
partial result in that direction is  Exercise 2(a) on page 292. 
 
Recall that the projective modules are precisely the direct summands 
of free modules.

By the classic theorem of Sabbagh-Eklof \cite{SaEk}, cf.\ 
Prest \cite[Theorem 14.25]{Pr:book} (which is stated there for 
the case of \emph{right} modules), \( R \) is left perfect  and
right coherent if{}f the class of projective left
\( R  \)-modules is elementary. The underlying algebraic result is 
due to Chase \cite{Ch}.

\begin{thm}
   Given\nlabel{Free:Mod:ai} a ring \(R\) of cardinality \(\tau\)\,,
   \(R^{(\tau^{+})}\) \tupop regarded as a left \(R\)-module\tupcp\ 
	is \(\tau^{+}\)-saturated if
   and only if \(R\) is left perfect and right coherent. 
\end{thm}

\begin{proof}
	 
	\noindent\( [\Rightarrow] \)\ Assume that 
	\(\lsub{R}{R^{(\tau^{+})}}\) is
	saturated (i.e. \(\tau^{+}\)-saturated).  By Theorem \ref{AlmInd:Ch} 
	every \(R^{(\lambda)}\) for \(\lambda \geq \tau^{+}\) is saturated (and
	elementarily equivalent to \(R^{(\tau^{+})}\)).  
	Let \(T = Th(R^{(\tau^{+})})\)\,. 
	As observed after Example \ref{Free:Mod} above \(T\) is tt. 

	We claim that the projective (left) \(R\)-modules are precisely the
	direct summands of models of \(T\)\,.  
	This will suffice, as
	by \cite[Lemma 2.23(a)]{Pr:book} this class coincides with the 
	class of pure submodules of models of \( T \)\,, and since 
	\( T \) is closed under products, by  \cite[Lemma 2.31]{Pr:book} 
	this latter class is elementary.
	
	If \(\lsub{R}{M}\) is projective then it is
	a direct summand of  \(R^{(\lambda)}\)  for some 
	\(\lambda \geq\tau^{+}\)\,, and the latter is a model of \(T\)\,. 
	
	On the other hand
	suppose that \(M\) is a direct summand of a model \(N\) of \(T\)\,. 
	Take \(\lambda \geq \tau^{+}\) sufficiently large so that 
   \(N\) is isomorphic to an
	elementary substructure of \( R^{(\lambda)} \)\,.  As \(T\) is
	tt, \(N\) is a direct summand of \(R^{(\lambda)} \) whereby \(M\) is also
	a direct summand of \(R^{(\lambda)} \)\,.  So \(M\) is projective.

	\noindent\( [\Leftarrow] \)\ 
	Assume that \( R \) is  left perfect and right coherent (and so the 
	projective left modules form an elementary class).
	By \cite[Cor.\ 14.22]{Pr:book}, 
	\(T= Th(\lsub{R}{R^{(\tau^{+})}})\) is tt. 
	As \(T\) is tt it has a saturated model \( M \) in power 
	\( \tau^{+} \)\,. 
	By freeness, there is a surjection 
	\( R^{(\tau^{+})}\twoheadrightarrow M \)\,.
	By the assumption \(M\) is projective as it is elementarily 
	equivalent to a free module; so it is isomorphic to a 
	summand, hence to an elementary submodel, of 
	\( R^{(\tau^{+})} \)\,.
	As \(T\) is tt and nonmultidimensional any elementary extension of 
	a \( \kappa \)-saturated model (\( \kappa\ge \card{T}^{+} \)) is 
	\( \kappa \)-saturated. Hence in particular, \( R^{(\tau^{+})} \) 
	is saturated. 
\end{proof}

\begin{remnc}
    Hence in particular if \( R \) is commutative and for some 
    \( \kappa>\card{R}+\aleph_{0} \)\,, \( R^{(\kappa)} \) is 
    saturated, then it has finite Morley  rank.
\end{remnc}

\begin{exm}
    By\nlabel{Free:ModUTMR}
    contrast, a well-known source of counter-examples
	 (cf. Small \cite{MR0188252}) provides us 
    with a saturated free module with infinite Morley rank, and so a 
    counter-example to the question of Pillay and Sklinos 
    \cite[Question 3.14]{MR3411166}. These examples are almost always 
	 presented for \emph{right} modules and we follow that custom here.
    
   Consider the upper triangular matrix ring 
    \[ 
    R=
    \left(
    \begin{array}
    {cc}
       
       \QQ & \QQ(x)  \\
       
       0 & \QQ(x)  \\
    \end{array}
    \right)
    \]
	 \( R \) is well-known to be right artinian and right perfect, 
	 but  only 
    left coherent,  not left noetherian. So by Theorem 
    \ref{Free:Mod:ai}, the free (right) module on \( \aleph_{1} \) 
	 generators is saturated.
    
    It is a standard exercise to determine \emph{all} the left and 
    right ideals of a ring of this sort.
    
    Furthermore, the pp-definable 
    subgroups of \( R_{R} \) are exactly the finitely generated 
    \emph{left} ideals of \( R \)\,, cf.\ 
    \ Prest \cite[Theorem 14.16]{Pr:book}: if \( I  \) is  
    generated by 
    \( \Set{r_{1},\,\ldots,\,r_{n}} \)\,, then it is defined in
    \( R_{R} \) by 
    \( \varphi(v)=\exists w_{1},\,\ldots,\,w_{n}\, 
    (v=w_{1}r_{1}+\cdots+w_{n}r_{n}) \)\,. 
	 
	 The Jacobson radical is the ``upper right corner'' of the ring; 
	 it is generated 
	 as a right ideal by any non-zero element, but as a left ideal it 
	 has the structure of \(\QQ(x) \) as a vector space over 
	 \( \QQ \)\,. In particular every finite dimensional 
	 \( \QQ \)-subspace \( V \) of \( \QQ(x) \) yields a pp-definable 
	 subgroup
	 \( \left(
	 \begin{array}
	 {cc}
	 	 
	 	 0 & V  \\
	 	 
	 	 0 & 0  
	 \end{array}\right)
	 \) of the  Jacobson radical, and hence we obtain infinite increasing 
	 chains of pp-definable subgroups.
	 
    Morley rank equals \( U \)-rank,  
    cf.\ Prest \cite[Theorem 5.18]{Pr:book}, and 
    the \( U \)-rank of a definable subgroup is just its rank in the 
    lattice of pp-definable subgroups,  cf.\ Prest
    \cite[Theorem 5.12]{Pr:book}, so the Morley rank of \( R_{R} \) 
	 is infinite. [In fact, a more careful analysis shows that the 
	 Morley rank is \( \omega+1 \)\,.]
    
	  %
\end{exm}
 

\section{Questions and open problems}

\begin{question}
    Is there a fundamental difference between theories that are
    \( (\tau,\tau^{+}) \) almost indiscernible and those that are
    \( (\mu,\tau^{+}) \) almost indiscernible 
    for some \( \mu<\tau \)\,?
\end{question}

In section \ref{Free:Th} we get a couple of results giving 
characterizations of algebraic closure (\ref{Free:310}) and 
independence in a saturated free algebra (\ref{Free:312}).
The next two questions relate to these results.

\begin{question}
    Are\nlabel{Q:FA-aind} there similar results for arbitrary almost 
    indiscernible theories?
\end{question}

\begin{question}
    Is\nlabel{Q:ind}
    there a more general description of independence in the theory 
    of a saturated free algebra?
    
    We are thinking of something that might fit into a general 
    abstract framework similar to that developed for theories of 
    modules in Prest \cite[\S5.4]{Pr:book}  (largely based on 
     \cite[Pillay-Prest]{PiPr83}).
    
    The closest direct analogue of Prest 
    \cite[Theorem 5.35]{Pr:book} 
    would be the following:
    
    \( \ola \) and \( \olb \) are independent over \( \olc \) if{}f 
    there is a basis \( X \)\,, the disjoint union of 
    \( A \)\,, \( B \)\,, and \( C \)\,, such that \( \olc\in \gen{C} \) 
    and \( \ola\in \gen{A\cup C} \)\,, 
    \( \olb\in \gen{B\cup C} \)\,.
\end{question}

\begin{question}
    Is\nlabel{Q:class} there any kind of classification of those 
    varieties \( \mathcal{V}\) for which
    the free algebra on \(\tau^{+}\) generators is 
    \(\tau^{+}\)-saturated?
    
    One should be cautious, as there are examples in Baldwin-Shelah 
    of such \( \mathcal{V}\)
    which have unstable algebras in the variety.
\end{question}

\begin{question}
    What\nlabel{Q:ClassSt}
    about Question \ref{Q:class}, assuming the stability of
    \(\mathcal{V}\)\,,  that is, that every
    completion of \(\Th(\mathcal{V})\) is stable?
\end{question}

\begin{question}
    In Propostion \ref{Free:Mod:max} we showed that the rank of the 
    type of a basic element of a large saturated free module is 
    maximal. Is this true for large  
    saturated free algebras in general?
\end{question}

For the rest, let us assume that the free algebra \(M\)  on
\(\tau^{+}\)-generators is \(\tau^{+}\) saturated, 
and let \(T = \Th(M)\)\,.

\begin{question}
   Is\nlabel{Q:tt} \( T \) totally transcendental? 
\end{question}

\begin{question}
    Is\nlabel{Q:struct} there a structure theorem
    for the algebra \( M \), 
    for example as
    some kind of a product of a module and of a combinatorial part, 
    along the lines of Hart and Valeriote
    \cite{MR1129148}?
\end{question}

\begin{question}
    Implicit\nlabel{Q:QE} in the last few questions is the following:
    
    Is there some kind of relative quantifier elimination theorem for 
    such theories?
\end{question}


\begin{center}
\begin{minipage}{5in}
\begin{center}ACKNOWLEDGEMENTS\end{center}
Work on the results reported here began when Thomas Kucera visited
Anand Pillay at Notre Dame during a sabbatical leave in late 2017, and
Dr.~Kucera would like to thank the University of Notre Dame for its
hospitality and support.  Dr.~Pillay would like to thank the NSF for
its support through grants DMS-1665035 and 1760413.  Dr.~Pillay would
also like to thank the University of Manchester for their support as a
Kathleen Ollerenshaw Professor in summer 2019.  Both authors thank
Mike Prest for very helpful discussions concerning the model theory of
modules, in particular the material related to Theorem 
\ref{Free:Mod:ai} and Example \ref{Free:ModUTMR}. The authors also 
want to express special gratitude to the two referees, who gave 
careful and insightful readings to several versions of the manuscript.
\end{minipage}
\end{center}



\end{document}